\documentclass[11pt,leqno]{amsart}
\pdfoutput=1
\usepackage{amsmath, amssymb}
\usepackage{graphicx,color,hyperref}
\usepackage{verbatim,enumitem,ifpdf}

\def\smallskip{\vskip4pt plus1pt minus1pt}
\def\medskip{\vskip7pt plus2pt minus1pt}
\def\bigskip{\vskip11pt plus2pt minus1pt}

\def\plainsubsection#1|{%
  \par\vskip0.25cm\penalty -100
  \noindent{\bf #1}
  \vskip3pt plus 1pt minus 0pt
  \penalty 500}

\long\def\claim#1|#2\endclaim{\par\vskip 4pt\noindent 
{\bf #1.}\ {\em #2}\par\vskip 5pt}

\def\today{\ifcase\month\or
January\or February\or March\or April\or May\or June\or July\or August\or
September\or October\or November\or December\fi \space\number\day,
\number\year}

\catcode`\@=11
\newcount\@tempcnta \newcount\@tempcntb 
\def\timeofday{{%
\@tempcnta=\time \divide\@tempcnta by 60 \@tempcntb=\@tempcnta
\multiply\@tempcntb by -60 \advance\@tempcntb by \time
\ifnum\@tempcntb > 9 \number\@tempcnta:\number\@tempcntb
\else\number\@tempcnta:0\number\@tempcntb\fi}}

\def\openup{\afterassignment\@penup\dimen@=}
\def\@penup{\advance\lineskip\dimen@
  \advance\baselineskip\dimen@
  \advance\lineskiplimit\dimen@}
\newdimen\jot \jot=3pt
\newskip\plaincentering \plaincentering=0pt plus 1000pt minus 1000pt
\def\ialign{\everycr{}\tabskip\z@skip\halign}
\def\plaineqalign#1{\null\,\vcenter{\openup\jot\m@th
  \ialign{\strut\hfil$\displaystyle{##}$&$\displaystyle{{}##}$\hfil
      \crcr#1\crcr}}\,}
\newif\ifdt@p
\def\displ@y{\global\dt@ptrue\openup\jot\m@th
  \everycr{\noalign{\ifdt@p \global\dt@pfalse \ifdim\prevdepth>-1000\p@
      \vskip-\lineskiplimit \vskip\normallineskiplimit \fi
      \else \penalty\interdisplaylinepenalty \fi}}}
\def\@lign{\tabskip\z@skip\everycr{}} 
\def\displaylines#1{\displ@y \tabskip\z@skip
  \halign{\hbox to\displaywidth{$\@lign\hfil\displaystyle##\hfil$}\crcr
    #1\crcr}}
\def\plaineqalignno#1{\displ@y \tabskip\plaincentering
  \halign to\displaywidth{\hfil$\@lign\displaystyle{##}$\tabskip\z@skip
    &$\@lign\displaystyle{{}##}$\hfil\tabskip\plaincentering
    &\llap{$\@lign##$}\tabskip\z@skip\crcr
    #1\crcr}}
\def\plainleqalignno#1{\displ@y \tabskip\plaincentering
  \halign to\displaywidth{\hfil$\@lign\displaystyle{##}$\tabskip\z@skip
    &$\@lign\displaystyle{{}##}$\hfil\tabskip\plaincentering
    &\kern-\displaywidth\rlap{$\@lign##$}\tabskip\displaywidth\crcr
    #1\crcr}}
\def\plaincases#1{\left\{\,\vcenter{\normalbaselines\m@th
    \ialign{$##\hfil$&\quad##\hfil\crcr#1\crcr}}\right.}
\def\plainmatrix#1{\null\,\vcenter{\normalbaselines\m@th
    \ialign{\hfil$##$\hfil&&\quad\hfil$##$\hfil\crcr
      \mathstrut\crcr\noalign{\kern-\baselineskip}
      #1\crcr\mathstrut\crcr\noalign{\kern-\baselineskip}}}\,}

\catcode`\@=12

-lmbx10 at 20.736pt
-lmbx10 at 14.4pt
\font\twelvebf = ec-lmbx10 at 12pt

\font\twelvebsy=cmbsy10 at 12pt
\font\tenbsy=cmbsy10
\font\eightbsy=cmbsy8
\font\sevenbsy=cmbsy7
\font\sixbsy=cmbsy6
\font\fivebsy=cmbsy5

\font\tenCal=eusm10 at 10pt
\font\eightCal=eusm10 at 8pt
\font\sevenCal=eusm10 at 7pt
\font\sixCal=eusm10 at 6pt
\font\fiveCal=eusm10 at 5pt

\newfam\Calfam
  \textfont\Calfam=\tenCal
  \scriptfont\Calfam=\sevenCal
  \scriptscriptfont\Calfam=\fiveCal
\def\Cal{\fam\Calfam\tenCal}

\font\tenmsa = msam10
\font\sevenmsa = msam7
\font\fivemsa = msam5
\newfam\msafam
\textfont\msafam=\tenmsa
\scriptfont\msafam=\sevenmsa
\scriptscriptfont\msafam=\fivemsa

\font\tenmsb = msbm10
\font\sevenmsb = msbm7
\font\fivemsb = msbm5
\newfam\msbfam
\textfont\msbfam=\tenmsb
\scriptfont\msbfam=\sevenmsb
\scriptscriptfont\msbfam=\fivemsb

\font\teneuf = eufm10
\font\seveneuf = eufm7
\font\fiveeuf = eufm5
\newfam\euffam
\textfont\euffam=\teneuf
\scriptfont\euffam=\seveneuf
\scriptscriptfont\euffam=\fiveeuf

\catcode`\@=11
\def\frak#1{{\fam\euffam#1}}
\def\Bbb#1{{\fam\msbfam#1}}

\def\hexnumber@#1{\ifcase#1 0\or1\or2\or3\or4\or5\or6\or7\or8\or9\or
 A\or B\or C\or D\or E\or F\fi}
\edef\msa@{\hexnumber@\msafam}
\edef\msb@{\hexnumber@\msbfam}
\mathchardef\square="0\msa@03
\mathchardef\smallsetminus="2\msb@72
\let\ssm=\smallsetminus
\mathchardef\compact="3\msa@62
\mathchardef\complement="3\msa@7B
\mathchardef\upharpoonright="3\msa@16

\mathchardef\subsetneq="3\msb@28
\mathchardef\supsetneq="3\msb@29
\catcode`\@=12

\def\ii{{\rm i}\kern0.8pt}
\def\bu{{\scriptstyle\bullet}}

\def\smallvee{{\scriptscriptstyle\vee}}

\def\build#1^#2_#3{\mathrel{\mathop{\null#1}\limits^{#2}_{#3}}}
\def\buildo#1\over#2{\mathrel{\mathop{\null#2}\limits^{#1}}}
\def\buildu#1\under#2{\mathrel{\mathop{\null#2}\limits_{#1}}}

\font\twelvemib = lmmib10 at 12pt
\font\tenmib = lmmib10
\font\sevenmib = lmmib10 at 7pt

\font\eightrm = ec-lmr10 at 8pt
\font\eightbf = ec-lmbx10 at 8pt
\font\eightit = ec-lmri10 at 8pt
\font\eighttt = ec-lmtt10 at 8pt
\font\eighti = lmmi10 at 8pt
\font\eightsy = lmsy10 at 8pt
\font\eightmib = lmmib10 at 8pt
\font\eightex = lmex10 at 8pt

\font\eightmsa = msam8
\font\eightmsb = msbm8
\font\eighteuf = eufm8

\font\sixrm = ec-lmr10 at 6pt
\font\sixbf = ec-lmbx10 at 6pt
-lmri10 at 6pt
-lmtt10 at 6pt
\font\sixi = lmmi10 at 6pt
\font\sixsy = lmsy10 at 6pt
\font\sixmib = lmmib10 at 6pt

\font\sixmsa = msam6
\font\sixmsb = msbm6
\font\sixeuf = eufm6

\font\fiverm = ec-lmr10 at 5pt
\font\fivebf = ec-lmbx10 at 5pt
-lmri10 at 5pt
-lmtt10 at 5pt
\font\fivei = lmmi10 at 5pt
\font\fivesy = lmsy10 at 5pt

\font\fivemib = lmmib10 at 5pt

\font\fivemsa = msam5
\font\fivemsb = msbm5

\def\eightpoint{\def\rm{\fam0\eightrm}
   \textfont0=\eightrm \scriptfont0=\sixrm \scriptscriptfont0=\fiverm
   \textfont1=\eighti  \scriptfont1=\sixi  \scriptscriptfont1=\fivei
   \textfont2=\eightsy \scriptfont2=\sixsy \scriptscriptfont2=\fivesy
   \textfont3=\eightex \scriptfont3=\eightex\scriptscriptfont3=\eightex
   \def\it{\fam\itfam\eightit}%
   \textfont\itfam=\eightit
   \def\bf{\fam\bffam\eightbf}%
   \textfont\bffam=\eightbf \scriptfont\bffam=\sixbf
   \scriptscriptfont\bffam=\fivebf
   \def\tt{\fam\ttfam\eighttt}%
   \textfont\ttfam=\eighttt
   \textfont\msafam=\eightmsa \scriptfont\msafam=\sixmsa
   \textfont\msbfam=\eightmsb \scriptfont\msbfam=\sixmsb
   \textfont\euffam=\eighteuf \scriptfont\euffam=\sixeuf
   \textfont\Calfam=\eightCal
   \scriptfont\Calfam=\sixCal
   \normalbaselineskip=9.6pt
   \normalbaselines\rm}

\newskip\eightpointsurround
\def\begineightpoint{\vskip\eightpointsurround \bgroup\eightpoint}
\def\endeightpoint{\vskip\eightpointsurround \egroup}

\def\twelvepointbf{%
 \textfont0=\twelvebf   \scriptfont0=\eightbf   \scriptscriptfont0=\sixbf
 \textfont1=\twelvemib  \scriptfont1=\eightmib \scriptscriptfont1=\sixmib
 \textfont2=\twelvebsy  \scriptfont2=\eightbsy  \scriptscriptfont2=\sixbsy
 \twelvebf
 \baselineskip=14.4pt}

\def\tenpointbf{%
 \textfont0=\tenbf\scriptfont0=\sevenbf\scriptscriptfont0=\fivebf
 \textfont1=\tenmib\scriptfont1=\sevenmib\scriptscriptfont1=\fivemib
 \textfont2=\tenbsy\scriptfont2=\sevenbsy\scriptscriptfont2=\fivebsy
 \tenbf\baselineskip=14.4pt}

\def\today{\ifcase\month\or
January\or February\or March\or April\or May\or June\or July\or August\or
September\or October\or November\or December\fi \space\number\day,
\number\year}

\catcode`\@=11
\newcount\@tempcnta \newcount\@tempcntb 
\def\timeofday{{%
\@tempcnta=\time \divide\@tempcnta by 60 \@tempcntb=\@tempcnta
\multiply\@tempcntb by -60 \advance\@tempcntb by \time
\ifnum\@tempcntb > 9 \number\@tempcnta:\number\@tempcntb
  \else\number\@tempcnta:0\number\@tempcntb\fi}}
\catcode`\@=12

\def\proof#1{\removelastskip\vskip4pt\noindent{\it #1.}}
\def\endproof{\vskip5pt plus 1pt minus 2pt}

\def\qed{\relax\ifmmode
            \hbox{$\square$}%
         \else
            {\unskip\nobreak\hfil
            \penalty50\hskip1em\null\nobreak\hfil\hbox{$\square$}%
            \parfillskip=0pt\finalhyphendemerits=0\endgraf}%
         \fi}
\def\bigsquare{{\kern-0.3ex\hbox{
\vrule height 1.7ex  width 0.093ex  depth 0ex\kern-0.093ex
\vrule height 1.8ex  width 1.7ex  depth -1.707ex\kern-0.093ex
\vrule height 1.7ex  width 0.093ex  depth 0ex\kern-1.65ex
\vrule height 0.093ex  width 1.6ex  depth 0ex}\kern0.3ex}}

\long\def\claim#1 #2\endclaim
{\removelastskip\vskip7pt plus 1pt minus 1pt\noindent{\bf#1.}
{\it\ignorespaces#2}\vskip\baselineskip}

\newdimen\plainitemindent\plainitemindent=18pt

\def\\{\hfil\break}

\def\bB{{\Bbb B}}
\def\bC{{\Bbb C}}
\def\bH{{\Bbb H}}
\def\bN{{\Bbb N}}

\def\bR{{\Bbb R}}
\def\bZ{{\Bbb Z}}

\def\cB{{\Cal B}}
\def\cC{{\Cal C}}
\def\cE{{\Cal E}}
\def\cF{{\Cal F}}
\def\cG{{\Cal G}}
\def\cH{{\Cal H}}
\def\cI{{\Cal I}}

\def\cL{{\Cal L}}
\def\cO{{\Cal O}}
\def\cU{{\Cal U}}
\def\cX{{\Cal X}}
\def\cXp{{\overline\cX{\kern0.75pt}'}}
\def\cY{{\Cal Y}}
\def\cZ{{\Cal Z}}

\def\gm{{\frak m}}

\def\stimes{\mathop{\kern0.7pt
\vrule height 0.4pt depth 0pt width 5pt\kern-5pt
\vrule height 5.4pt depth -5pt width 5pt\kern-5pt
\vrule height 5.4pt depth 0pt width 0.4pt\kern4.6pt
\vrule height 5.4pt depth 0pt width 0.4pt\kern-6.5pt
\raise0.3pt\hbox{$\times$}\kern-0.7pt}\nolimits}

\catcode`@=11
\def\rightarrowfil{\m@th\mathord-\mkern-6mu%
  \cleaders\hbox{$\mkern-2mu\mathord-\mkern-2mu$}\hfill
  \mkern-6mu\mathord\rightarrow}
\catcode`@=12

\def\dbar{\overline\partial}
\def\ddbar{\partial\overline\partial}
\def\Ker{\mathop{\rm Ker}\nolimits}
\def\End{\mathop{\rm End}\nolimits}
\def\GL{\mathop{\rm GL}\nolimits}
\def\Herm{\mathop{\rm Herm}\nolimits}
\def\Id{\mathop{\rm Id}\nolimits}

\def\Re{\mathop{\rm Re}\nolimits}

\def\mod{\mathop{\rm mod}\nolimits}
\def\diam{\mathop{\rm diam}\nolimits}
\def\codim{\mathop{\rm codim}\nolimits}

\def\Aut{\mathop{\rm Aut}\nolimits}
\def\Hom{\mathop{\rm Hom}\nolimits}

\def\plainitem#1{\par\parindent=\plainitemindent\noindent
\hangindent\parindent\hbox to\parindent{#1\hss}\ignorespaces}

\def\plainsection#1#2{%
  \par\vskip0.75cm\penalty -100
  \vbox{\baselineskip=14.4pt\noindent{{\twelvepointbf #1.\kern6pt #2}}}
  \vskip5pt\penalty 500}

\def\subsection#1#2{%
  \par\vskip0.2cm\penalty-100
  \vbox{\noindent{{\tenpointbf #1. #2}}}
  \vskip1pt
  \penalty 500}

\def\subsubsection#1#2{%
  \par\vskip0.1cm\penalty -100
  \noindent{{\tenpointbf #1.} {\it #2}}}

\def\Bibitem#1&#2&#3;&#4&%
{\hangindent=1.75cm\hangafter=1
\noindent\rlap{\hbox{\bf #1}}\kern1.75cm{\rm #2}{\it #3}\/; {\rm #4.}}

\def\setheadtitles#1{%
\setbox\sectionbox\hbox{\eightpoint #1}
\setbox\titlebox\hbox{\eightpoint #1}}

\let\em=\it 
\long\def\specnote#1#2{#1\footnote{}{\baselineskip=9.6pt{%
\kern-\parindent\eightrm #1 #2\par\vskip-\baselineskip}}}

\long\def\claim#1|#2\endclaim{\par\vskip 6pt\noindent 
{\bf #1.}\ {\em #2}\par\vskip 5pt}

\def\qed{\hfill\hbox{
\vrule height 1.453ex  width 0.093ex  depth 0ex
\vrule height 1.5ex  width 1.3ex  depth -1.407ex\kern-0.1ex
\vrule height 1.453ex  width 0.093ex  depth 0ex\kern-1.35ex
\vrule height 0.093ex  width 1.3ex  depth 0ex}}

\def\lreqno#1#2{\leqno\hbox to 0.01pt{\rlap{\hbox{$#1$}}}\kern-0.02pt
\rlap{\kern\hsize\kern 0.01pt\hbox to 0.01pt{\llap{\hbox{$#2$}}}}}

\let\ol=\overline

\let\wt=\widetilde
\def\loc{{\rm loc}}

\def\pr{\mathop{\rm pr}\nolimits}
\def\Id{\mathop{\rm Id}\nolimits}

\def\exph{\mathop{\rm exph}\nolimits}
\def\logh{\mathop{\rm logh}\nolimits}

\title[Bergman bundles and applications]{Bergman bundles and applications
to the\vskip4pt geometry of compact complex manifolds}

\author[Jean-Pierre Demailly]{Jean-Pierre Demailly}

\begin{document}

\begin{abstract} We introduce the concept of Bergman bundle attached
to a hermitian manifold $X$, assuming the manifold $X$ to be compact --
although the results are local for a large part. The Bergman bundle
is some sort of infinite dimensional very ample Hilbert bundle
whose fibers are isomorphic to the standard $L^2$ Hardy space on the complex
unit ball; however the bundle is locally trivial only in the real analytic
category, and its complex structure is strongly twisted. We compute
the Chern curvature of the Bergman bundle, and show that it is strictly
positive. As a potential application, we investigate a long standing
and still unsolved conjecture of Siu on the invariance of plurigenera
in the general situation of polarized families of compact K\"ahler manifolds.

\medskip

\vskip3pt\noindent
{\sc MSC Classification 2020.} 32J25, 32F32

\vskip3pt\noindent
{\sc Funding.} The author is supported by the Advanced ERC grant
ALKAGE, no 670846 from September 2015, attributed by the
European Research Council.
\end{abstract}

\maketitle

\hbox to \textwidth{\hfill\it Dedicated to Professor Bernard Shiffman
on the occasion of his retirement}

\plainsection{0}{Introduction}

Projective varieties are characterized, almost by definition, by the existence
of an ample line bundle. By the Kodaira embedding theorem [Kod54], they are
also characterized among compact complex manifolds by the existence of a
positively curved holomorphic line bundle, or equivalently, of a Hodge metric,
namely a K\"ahler metric with rational cohomology class. On the other hand,
general compact K\"ahler manifolds, and especially general complex tori, fail to
have a positive line bundle. Still, compact K\"ahler manifolds
possess topological complex line bundles of positive curvature, that are
in some sense arbitrary close to being holomorphic, see e.g.\
[Lae02] and [Pop13]. It may nevertheless come as a surprise that every compact
complex manifold carries some sort of very ample holomorphic vector
bundle, at least if~one accepts certain Hilbert bundles of infinite
dimension. Motivated by geometric quantization, Lempert and Sz\H{o}ke
[LeS14] have introduced and discussed a more general concept of
``field of Hilbert spaces'' which is similar in spirit.

\claim 0.1. Theorem|Every compact complex manifold $X$ carries a locally
trivial real analytic Hilbert bundle $B_\varepsilon\to X$ of
infinite dimension, defined for $0<\varepsilon
\le\varepsilon_0$, equipped with an integrable $(0,1)$-connection
$\overline\partial=\nabla^{0,1}$ $($in a generalized sense$)$, that is
a closed densely defined operator in the space of $L^2$ sections,
in such a way that the sheaf $\cB_\varepsilon=\cO_{L^2}(B_\varepsilon)$ of
$\overline\partial$-closed locally $L^2$ sections is ``very ample'' in the
following sense.
\plainitem{\rm(a)} $H^q(X,\cB_\varepsilon\otimes_\cO\cF)=0$ for every
$(\,$finite rank$\,)$ coherent sheaf $\cF$ on $X$ and every $q\ge 1.$
\plainitem{\rm(b)} Global sections of the Hilbert space
$\bH=H^0(X,\cB_\varepsilon)$ provide an embedding of $X$ into a certain
Grassmannian of closed subspaces of infinite codimension in $\bH$.
\plainitem{\rm(c)}
The bundle $B_\varepsilon$ carries a natural Hilbert metric $h$
such that
the curvature tensor $\ii\Theta_{B_\varepsilon,h}$ is Nakano positive $($and
even Nakano positive unbounded$\;!)$.\vskip0pt
\endclaim

\noindent Parts (a) and (b) are proved by considering the (pre)sheaf
structure of $\cB_\varepsilon$, and observing that there is a related
$L^2$ Dolbeault complex on which H\"ormander's $L^2$ estimates [Hor66]
can be applied. The case of Stein manifolds is sufficient, and the
corresponding Hilbert bundle $B_\varepsilon$ is not involved in the
arguments. Technically, the proof is given in Proposition~2.5 and
Remark~2.6. Part~(c) deals with the geometry of $B_\varepsilon$, and is
treated in section~3.

We start by explaining a little bit more the relationship between the
``Hilbert bundles'' involved here, and the more familiar concept of
locally trivial holomorphic Hilbert bundle~: such a bundle $E\to X$ is
required to be trivial on sufficiently small open sets $V\subset X$, and
such that $E_{|V}\simeq V\times \cH$
where $\cH$ is a complex Hilbert space. The gluing transition
automorphism with another local trivialization
$E_{|V'}\simeq V'\times\cH$ should then be of the form $(z,\xi)\mapsto
(z,g(z)\cdot\xi)$ where $g$ is a holomorphic map from $V\cap V'$
to the open set $\GL(\cH)$ of invertible continuous operators in $\End(\cH)$.
Smooth and real analytic locally trivial Hilbert bundles can be defined
in a similar manner by requiring $g$ to be in $C^\infty$, resp.\ in $C^\omega$.
A~smooth hermitian structure on $E$ is a smooth family of hermitian
metrics $h(z)$ on the fibers, given in the trivializations by smooth maps
\hbox{$h_V\in C^\infty(V,\Herm_+(\cH))$}, where $\Herm_+(\cH)$ is the set of
positive definite (coercive) hermitian forms on $\cH$. The usual formalism
of Chern connections still applies: one gets a unique
connection $\nabla_h=\nabla_h^{1,0}+\nabla_h^{0,1}$ acting on
$C^\infty(X,E)$ in such a way that $h$ is $\nabla_h$ parallel and 
$\nabla_h^{0,1}=\overline\partial$; moreover, the kernel of $\nabla_h^{0,1}$
coincides with the sheaf of holomorphic sections $\cO_X(E)$. This connection
is given by exactly the same formulas as in the finite dimensional case,
namely \hbox{$\nabla_h^{1,0}\simeq h_V^{-1}\circ \partial\circ h_V$}
over $V$, with a curvature tensor $\nabla_h^2=\Theta_{E,h}$ given by
$\Theta_{E,h}\simeq\overline\partial (h_V^{-1}\partial h_V)$
(if one views $h_V(z)$ as an endomorphism of $\cH$); locally, $\Theta_{E,h}$
can thus be seen as a smooth $(1,1)$-form with values in the space of
continuous endomorphisms $\End(\cH)$. In general, if $E$ is a
smooth~Hilbert bundle (defined as above, but with gluing automorphisms
$g\in C^\infty(V,\GL(\cH))$), a smooth $(0,1)$-connection $\nabla_A^{0,1}$ is
an order $1$ linear differential operator that is locally of the
form $\overline\partial+A_V$ where
$A_V\in C^\infty(V,\Lambda^{0,1}T^*_X\otimes\End(\cH))$. It is said to
be integrable if $(\nabla_A^{0,1})^2=0$, i.e.\
$\overline\partial A_V+A_V\wedge A_V=0$ on each trivializing chart~$V$.
The~following equivalence of categories is well known, and follows e.g.\
from Malgrange [Mal58, chap.~X, Theorem 1], although the statement is
expressed there in more concrete terms.

\claim 0.2. Theorem {\rm (Malgrange [Mal58])}|The category of holomorphic
vector bundles on $X$ is equivalent to the category of smooth bundles
equipped with smooth integrable $(0,1)$-connec\-tions $\nabla_A^{0,1}$,
the holomorphic structure being obtained by taking the kernel sheaf
of~$\nabla_A^{0,1}$.
\endclaim

\noindent
Notice that the usual finite dimensional proofs apply essentially unchanged
to the case of locally trivial Hilbert bundles. For instance, one can adapt
Malgrange's inductive proof [Mal58] based on the Cauchy formula in one
variable (for holomorphic functions with values in a Banach space, depending
smoothly on some other parameters), or use a Nash-Moser process along
with the Bochner-Martinelli kernel (see e.g.\ [Web89]), or an
infinite dimensional version of H\"ormander's $L^2$ estimates (the latter
do not depend on the rank of bundles and are thus valid for Hilbert
bundles equipped with integrable smooth $(0,1)$-connections; the solution
of minimal $L^2$ norm can be used to find local $\nabla_A^{0,1}$-closed sections
generating fibers of the bundle). The result is also
valid for the category of real analytic Hilbert bundles, assuming $E$ and
$\nabla_A^{0,1}$ to be real analytic; the resulting holomorphic structure
on $E$ is then compatible with the originally given real analytic structure.
Now, when $X$ is compact hermitian, the standard $L^2$ techniques of
PDE theory lead to considering the space $L^2(X,E)$ of $L^2$ sections.
A~smooth $(0,1)$-connection then gives rise to a closed densely defined
operator
$$
\nabla_A^{0,1}:L^2(X,E)\longrightarrow L^2(X,\Lambda^{0,1}T^*_X\otimes E)
\leqno(0.3)
$$
that is never continuous. In our situation, the bundles come in a natural
way as a family of smooth (and even real analytic) Hilbert bundles
$E_\varepsilon$, $0<\varepsilon\leq\varepsilon_0$, associated
with a family $\cH_\varepsilon$ of Hilbert spaces which form a ``scale'',
in the sense that there are continuous injections with dense image
$\cH_{\varepsilon'}\hookrightarrow\cH_\varepsilon$,
$0<\varepsilon<\varepsilon'\leq\varepsilon_0$. The transition
automorphisms defining the $E_\varepsilon$'s are supposed to come from
invertible automorphisms of $\cH_{>0}=\bigcup_{\varepsilon>0}\cH_\varepsilon$
preserving each $\cH_\varepsilon$
(here $\cH_{>0}$ is just an inductive limit of Hilbert spaces).
Then it makes sense to consider generalized $(0,1)$-connections
that are locally of the form
$$
\nabla^{0,1}_A\simeq\overline\partial + A_V,\qquad
A_V\in C^\infty(\Lambda^{0,1}T^*_X\otimes\End(\cH_{>0})),
\leqno(0.4)
$$
where we actually have $A_{V|\cH_{\varepsilon'}}
\in C^\infty(\Lambda^{0,1}T^*_X\otimes
\Hom(\cH_{\varepsilon'},\cH_\varepsilon))$ for all
$0<\varepsilon<\varepsilon'\leq\varepsilon_0$. By our assumptions,
such connections still induce densely defined operators on each of the
spaces $L^2(X,E_\varepsilon)$, and we declare them to be integrable
when $(\nabla^{0,1}_A)^2=0$. The usual algebraic formalism for extending
the connection to higher degree forms and calculating the curvature tensor
still applies in this setting.

However, it may happen, and this will be the case for the Chern connection
matrices of our bundles $B_\varepsilon$ of Theorem 0.1,
that the $A_V$ do not induce continuous endomorphisms of
$\cH_\varepsilon$ (for any value of $\varepsilon>0$), although the kernel
of $\nabla^{0,1}$ in $L^2(X,B_\varepsilon)$ looks very much like a space
of holomorphic sections. In this context, the associated curvature
tensor $\Theta_{E_\varepsilon,h}$ need not either
take values in the continous endomorphisms. Then Malgrange's theorem implies
that such bundles do not correspond to locally trivial holomorphic bundles
as defined above, even under the integrability assumption.
At the end of Section 3 we will briefly discuss in which sense
$B_\varepsilon$ can still be considered to be some sort of infinite
dimensional complex space, in a way that the projection
map $B_\varepsilon\to X$ becomes holomorphic.

The cons\-truction of $B_\varepsilon$ is made by embedding $X$
diagonally in $X\times\ol X$ and taking a Stein tubular neighborhood
$U_\varepsilon$ of the diagonal, according to a well known technique
of Grauert [Gra58]. When $U_\varepsilon$ is chosen to be a geodesic
neighborhood with respect to some real analytic hermitian metric,
one can arrange that the first projection $p:U_\varepsilon\to X$
is a real analytic bundle whose fibers are biholomorphic to
hermitian balls. One then takes $B_\varepsilon$ to be a ``Bergman bundle'',
consisting of holomorphic $n$-forms $f(z,w)\,dw_1\wedge\ldots\wedge dw_n$
that are $L^2$ on the fibers $p^{-1}(z)\simeq B(0,\varepsilon)$. The fact that
$U_\varepsilon$ is Stein and real analytically locally trivial over $X$
then implies Theorem~0.1, using the corresponding Bergman type Dolbeault
complex.

In [Ber09], given a holomorphic fibration $\pi:X\to Y$
and a positive hermitian holomorphic line bundle $L\to X$,
Berndtsson has introduced a formally similar $L^2$ bundle $Y\ni t\mapsto A^2_t$,
whose fibers consist of sections of the adjoint bundle $K_{X/Y}\otimes L$
on the fibers $X_t=\pi^{-1}(t)$ of $\pi$, equipped with the corresponding
Bergman metric. In the situation considered by Berndtsson, the major
application is the case when $\pi$ is proper, so that $A^2_t$ is
finite dimensional, and is the holomorphic bundle associated with
the direct image sheaf $\pi_*(\cO_X(K_{X/Y}\otimes L))$.
The main result of [Ber09] is a calculation of the curvature,
and a proof that the direct image is a Nakano positive vector bundle.
On the other hand, when $\pi:X\to Y$ is non proper, and
especially when $(X_t)$ is a smooth family of smoothly bounded Stein
domains, the corresponding spaces $A^2_t$ are infinite dimensional
Hilbert spaces. The curvature of the corresponding Hilbert bundle has
been obtained by Wang Xu [Wan17] in this general setting.
Our curvature calculations can be seen as the very special case
where the fibers are smoothly varying hermitian balls and the centers
vary antiholomorphically. The calculation
can then be made in a very explicit way, by first considering the model case
of balls of constant radius in $\bC^n$, and then by using an osculation and
suitable Taylor expansions, in the case of varying hermitian metrics
(a similar osculating technique has been used in [ZeZ18] for the
study of Bargmann-Fock spaces). As a consequence, we get

\claim 0.5. Proposition|The curvature tensor of $(B_\varepsilon,h)$ admits 
an asymptotic expansion
$$
\langle(\Theta_{B_\varepsilon,h}\,\xi)(v,Jv),\xi\rangle_h
=\sum_{p=0}^{+\infty}\varepsilon^{-2+p}Q_p(z,\xi\otimes v),
$$
where, in suitable normal coordinates, the leading term
$Q_0(z,\xi\otimes v)$ is exactly equal to the curvature tensor of
the Bergman bundle associated
with the translation invariant tubular neighborhood
$$
U_\varepsilon=\{(z,w)\in\bC^n\times\bC^n\,;\;|z-\ol w|<\varepsilon\},
$$
in the ``model case'' $X=\bC^n$. That term $Q_0$ is an unbounded
quadratic hermitian form.
\endclaim

\noindent The potential geometric applications we have in mind are for
instance the study of Siu's conjecture on the K\"ahler invariance of
plurigenera (see 4.1 below), where the algebraic proof ([Siu02], [Pau07])
uses an auxiliary ample line bundle $A$. In the K\"ahler case at least,
one possible idea would be to replace $A$ by the infinite dimensional Bergman
bundle $B_\varepsilon$. The proof works to some extent, but some crucial
additional estimates seem to be missing to get the conclusion, see \S4.
Another question where Bergman bundles could potentially be useful
is the conjecture on transcendental Morse inequalities
for real $(1,1)$-cohomology classes $\alpha$ in the Bott-Chern cohomology
group $H^{1,1}_{\rm BC}(X,\bC)$. In that situation, multiples $k\alpha$ can be
approximated by a sequence of integral classes $\alpha_k$ corresponding
to topological line bundles $L_k\to X$ that are closer and closer to being
holomorphic, see e.g.\ [Lae02]. However, on the Stein tubular
neighborhood~$U_\varepsilon$, the pull-back
$p^*L_k$ can be given a structure of a genuine holomorphic line bundle
with curvature form very close to $k\,p^*\alpha$.
Our hope is that an appropriate Bergman theory of ``Hilbert dimension''
(say, in the spirit of Atiyah's $L^2$ index theory) can be used to
recover the expected Morse inequalities. There seem to be still
considerable difficulties in this direction, and we wish to leave
this question for future research.

The author addresses warm thanks to the referees for a number of useful
suggestions and observations that led to substantial improvements of the
original presentation.

\plainsection{1}{Exponential map and tubular neighborhoods}

Let $X$ be a compact $n$-dimensional complex manifold and $Y\subset X$ a
smooth totally real submanifold, i.e.\ such that $T_Y\cap JT_Y=\{0\}$ for
the complex structure $J$ on~$X$. By a well known result of Grauert [Gra58],
such a $Y$ always admits a fundamental system of Stein tubular neighborhoods
$U\subset X$ (this would be even true when $X$ is noncompact, but we only
need the compact case here). In fact, if $(\Omega_\alpha)$ is a finite
covering of $X$ such that $Y\cap\Omega_\alpha$ is a smooth complete
intersection $\{z\in\Omega_\alpha\,;\;x_{\alpha,j}(z)=0\}$, $1\leq j\leq q$
(where $q=\codim_\bR Y\geq n$), then
one can take $U=U_\varepsilon=\{\varphi(z)<\varepsilon\}$ where
$$
\varphi(z)=\sum_\alpha \theta_\alpha(z)\sum_{1\leq j\leq q}(x_{\alpha,j}(z))^2
\geq 0
\leqno(1.1)
$$
where $(\theta_\alpha)$ is a partition of unity subordinate to
$(\Omega_\alpha)$. The reason is that $\varphi$ is strictly plurisubharmonic
near~$Y$, as
$$
\ii\ddbar\varphi_{|Y}=
2i\sum_\alpha \theta_\alpha(z)\sum_{1\leq j\leq q}\partial x_{\alpha,j}
\wedge \dbar x_{\alpha,j}
$$
and $(\partial x_{\alpha,j})_j$ has rank $n$ at every point of $Y$, by the
assumption that $Y$ is totally real.

Now, let $\overline X$ be the complex conjugate manifold associated with
the integrable almost complex structure $(X,-J)$ (in other words,
$\cO_{\overline X}=\overline{\cO_X}$); we denote by $x\mapsto\overline x$
the identity map $\Id:X\to\overline X$ to stress that it is
conjugate holomorphic. The underlying real analytic manifold
$X^\bR$ can be embedded diagonally in $X\times\overline X$ by the
diagonal map $\delta:x\mapsto(x,\overline x)$, and the image
$\delta(X^\bR)$ is a totally real submanifold of $X\times \overline X$.
In fact, if $(z_{\alpha,j})_{1\leq j\leq n}$ is a holomorphic coordinate
system relative to a finite open covering $(\Omega_\alpha)$ of $X$,
then the $\overline z_{\alpha,j}$ define holomorphic coordinates
on $\overline X$ relative to $\overline\Omega_\alpha$, and
the ``diagonal'' $\delta(X^\bR)$ is the totally real submanifold
of pairs $(z,w)$ such that $w_{\alpha,j}=\overline z_{\alpha,j}$ for
all $\alpha,j$. In that case, we can take Stein tubular neighborhoods
of the form $U_\varepsilon=\{\varphi<\varepsilon\}$ where
$$
\varphi(z,w)=\sum_\alpha \theta_\alpha(z)\theta_\alpha(w)
\sum_{1\leq j\leq q}|\overline w_{\alpha,j}-z_{\alpha,j}|^2.
\leqno(1.2)
$$
Here, the strict plurisubharmonicity of $\varphi$ near $\delta(X^\bR)$ is
obvious from the fact that
$$
|w_{\alpha,j}-\overline z_{\alpha,j}|^2=
|z_{\alpha,j}|^2+|w_{\alpha,j}|^2-2\Re(z_{\alpha,j}w_{\alpha,j}).
$$
For $\varepsilon>0$ small, the first projection $\pr_1:U_\varepsilon\to X$
gives a complex fibration whose fibers are $C^\infty$-diffeomorphic to
balls, but they need not be biholomorphic to complex balls in general.
In order to achieve this property, we proceed in the following way. Pick a
real analytic hermitian metric $\gamma$ on $X\,$; take e.g.\ the $(1,1)$-part
$\gamma=g^{(1,1)}={1\over 2}(g+J^*g)$ of the Riemannian metric obtained
as the pull-back $g=\delta^*(\sum_j i df_j\wedge d\overline f_j)$,
where the $(f_j)_{1\leq j\leq N}$ provide a holomorphic immersion of the
Stein neighborhood $U_\varepsilon$ into~$\bC^N$. Let $\exp:T_X\to X$,
$(z,\xi)\mapsto \exp_z(\xi)$ be the exponential map associated with
the metric $\gamma$, in such a way that $\bR\ni t\mapsto\exp_z(t\xi)$ are
geodesics ${D\over dt}({du\over dt})=0$ for the
the Chern connection $D$ on $T_X$ (see e.g.\ [Dem94, (2.6)]).
Then $\exp$ is real analytic, and we have Taylor expansions
$$
\exp_z(\xi)=\sum_{\alpha,\beta\in\bN^n}a_{\alpha\beta}(z)\xi^\alpha
\smash{\overline\xi}^\beta,\qquad\xi\in T_{X,z}
$$
with real analytic coefficients $a_{\alpha\beta}$, where
$\exp_z(\xi)=z+\xi+O(|\xi|^2)$ in
local coordinates. The real analyticity means that these expansions are
convergent on a neighborhood $|\xi|_\gamma<\varepsilon_0$ of the zero
section of~$T_X$. We define the fiber-holomorphic part of the exponential
map to be
$$
\exph:T_X\to X,\qquad(z,\xi)\mapsto\exph_z(\xi)=
\sum_{\alpha\in\bN^n}a_{\alpha 0}(z)\xi^\alpha.\leqno(1.3)
$$
It is uniquely defined, is convergent on the same tubular neighborhood
$\{|\xi|_\gamma<\varepsilon_0\}$, has the property that
$\xi\mapsto\exph_z(\xi)$ is holomorphic for $z\in X$ fixed, and
satisfies again
$\exph_z(\xi)=z+\xi+O(\xi^2)$ in coordinates. By the implicit function,
theorem, the map $(z,\xi)\mapsto (z,\exph_z(\xi))$ is a real analytic
diffeomorphism from a neighborhood of the zero section of $T_X$ onto
a neighborhood $V$ of
the diagonal in $X\times X$. Therefore, we get an inverse real
analytic mapping $X\times X\supset V\to T_X$, which we denote by
$(z,w)\mapsto (z,\xi)$, $\xi=\logh_z(w)$, such that
$w\mapsto\logh_z(w)$ is holomorphic on $V\cap(\{z\}\times X)$, and
$\logh_z(w)=w-z+O((w-z)^2)$ in coordinates. The tubular neighborhood
$$
U_{\gamma,\varepsilon}=\{(z,w)\in X\times \overline X\,;~
|\logh_z(\overline w)|_\gamma<\varepsilon\}
$$
is Stein for $\varepsilon>0$ small; in fact, if $p\in X$ and
$(z_1,\ldots,z_n)$ is a holomorphic coordinate system centered at $p$
such that $\gamma_p=i\sum dz_j\wedge d\overline z_j$, then
$|\logh_z(\overline w)|_\gamma^2=|\overline w-z|^2+O(|\overline w-z|^3)$, hence
$\ii\ddbar|\logh_z(\overline w)|_\gamma^2>0$ at $(p,\overline p)\in
X\times\overline X$. By construction, the fiber $\pr_1^{-1}(z)$ of
$\pr_1:U_{\gamma,\varepsilon}\to X$ is biholomorphic to the
$\varepsilon$-ball of the complex vector space $T_{X,z}$ equipped with
the hermitian metric $\gamma_z$. In this way, we get a locally
trivial real analytic bundle $\pr_1:U_{\gamma,\varepsilon}$ whose fibers
are complex balls; it is important to notice, however, that this ball
bundle need not -- and in fact, will never -- be holomorphically
locally trivial.

\plainsection{2}{Bergman bundles and Bergman Dolbeault complex}

Let $X$ be a $n$-dimensional compact complex manifold equipped with a
real analytic hermitian metric~$\gamma$,
$U_\varepsilon=U_{\gamma,\varepsilon}\subset X\times\overline X$
the ball bundle considered in \S1 and
$$
p=(\pr_1)_{|U_\varepsilon}:U_\varepsilon\to X,\qquad
\overline p=(\pr_2)_{|U_\varepsilon}:U_\varepsilon\to\overline X
$$
the natural projections. We introduce what we call the ``Bergman
direct image sheaf''
$$
\cB_\varepsilon=p^{L^2}_*(\overline p^*\cO(K_{\overline X})).
\leqno(2.1)
$$
By definition, its space of sections $\cB_\varepsilon(V)$ over an
open subset $V\subset X$ consists of holomorphic sections $f$ of
$\overline p^*\cO(K_{\overline X})$ on $p^{-1}(V)$ that are
in $L^2(p^{-1}(K))$  for all compact subsets $K\compact V$,
i.e.
$$
\int_{p^{-1}(K)}i^{n^2}f\wedge \overline f\wedge \gamma^n
<+\infty,\quad\forall K\compact V.
\leqno(2.2)
$$
Then $\cB_\varepsilon$ is clearly a sheaf of infinite dimensional
Fr\'echet $\cO_X$-modules. In the case of finitely generated sheaves
over $\cO_X$,
there is a well known equivalence of categories between holomorphic vector
bundles $G$ over $X$ and locally free $\cO_X$-modules~$\cG$. As is well
known, the correspondence is given by $G\mapsto \cG:=\cO_X(G)={}$sheaf
of germs of holomorphic sections of~$G$, and the converse functor is
$\cG\mapsto G$, where $G$ is the holomorphic vector bundle whose fibers
are $G_z=\cG_z/\gm_z\cG_z=\cG_z\otimes_{\cO_{X,z}}\cO_{X,z}/\gm_z$ where
$\gm_z\subset\cO_{X,z}$ is the maximal ideal. In the case of
$\cB_\varepsilon$, we cannot take exactly the same route, mostly because the
desired ``holomorphic Hilbert bundle'' $B_\varepsilon$ will not even be
locally trivial in the complex analytic sense. Instead, we define directly
the fibers $B_{\varepsilon,z}$ as the set of holomorphic sections $f$
of $K_{\ol X}$ on the fibers $U_{\varepsilon,z}=p^{-1}(z)$, such that
$$
\int_{U_{\varepsilon,z}}i^{n^2}f\wedge \overline f<+\infty.
\leqno(2.2_z)
$$
Since $U_{\varepsilon,z}$ is biholomorphic to the unit ball $\bB_n\subset\bC^n$,
the fiber $B_{\varepsilon,z}$ is isomorphic to the Hilbert space
$\cH^2(\bB_n)$ of $L^2$ holomorphic $n$-forms on~$\bB_n$.
In fact, if we use orthonormal coordinates $(w_1,\ldots,w_n)$
provided by $\exph$ acting on the hermitian space
$(T_{X,z},\gamma_z)$ and centered at $\overline z$, we get
a biholomorphism $\bB_n\to p^{-1}(z)$ given by the
homothety $\eta_\varepsilon:w\mapsto \varepsilon w$,
and a corresponding isomorphism
$$
\plainleqalignno{
&B_{\varepsilon,z}\longrightarrow \cH^2(\bB_n),\qquad f\longmapsto
g=\eta_\varepsilon^*f,~~\hbox{i.e.\ with $I=\{1,\ldots,n\}$},&(2.3)\cr
\noalign{\vskip5pt}
&f_I(w)\,dw_1\wedge\ldots\wedge dw_n\longmapsto
\varepsilon^n\,f_I(\varepsilon w)\,dw_1\wedge\ldots\wedge dw_n,\quad
w\in\bB_n,&(2.3')\cr
&\Vert g\Vert^2=\int_{\bB_n}2^{-n}i^{n^2}g\wedge\overline g,\quad
g=g(w)\,dw_1\wedge\ldots\wedge dw_n\in\cH^2(\bB_n).&(2.3'')\cr}
$$
As $U_\varepsilon\to X$ is real analytically locally trivial over $X$,
it follows immediately that \hbox{$B_\varepsilon\to X$} is also a locally
trivial real analytic Hilbert bundle of typical fiber $\cH^2(\bB_n)$,
with the natural Hilbert metric obtained by declaring (2.3)
to be an isometry. Since $\Aut(\bB_n)$ is a real Lie group, the gauge group
of $B_\varepsilon\to X$ can be reduced to real analytic sections
of  $\Aut(\bB_n)$ and we have a well defined class of real analytic
connections on $B_\varepsilon$. In this context, one should
pay attention to the fact that a section $f$ in $\cB_\varepsilon(V)$
does not necessarily restrict to $L^2$ holomorphic sections
$f_{U_{\varepsilon,z}}\in B_{\varepsilon,z}$ for all $z\in V$, although
this is certainly true for almost all $z\in V$ by the Fubini theorem;
this phenomenon can already be seen through the fact that one does
not have a continuous
restriction morphism \hbox{$\rho_n:\cH^2(\bB_n)\to \cH^2(\bB_{n-1})$}
to the hyperplane
$z_n=0$. In fact, the function $(1-z_1)^{-\alpha}$ is
in $\cH^2(\bB_n)$ if and only
if $\alpha<(n+1)/2$, so that $(1-z_1)^{-n/2}$ is outside of the domain
of~$\rho_n$. As a consequence, the morphism $\cB_{\varepsilon,z}\to
B_{\varepsilon,z}$ (stalk of sheaf to vector bundle fiber) only has a
dense domain of definition, containing e.g.\ $\cB_{\varepsilon',z}$
for any $\varepsilon'>\varepsilon$. This is a familiar
situation in Von Neumann's theory of operators.

We now introduce a natural ``Bergman version'' of the Dolbeault
complex, by introducing a sheaf $\cF_\varepsilon^q$ over $X$ of
$(n,q)$-forms which can be written locally over small
open sets $V\subset X$ as
$$
f(z,w)=\sum_{|J|=q}f_J(z,w)\,dw_1\wedge\ldots\wedge dw_n\wedge
d\overline z_J,\quad
(z,w)\in U_\varepsilon\cap(V\times\overline X),
\leqno(2.4)
$$
where the $f_J(z,w)$ are $L^2_\loc$ smooth functions on
$U_\varepsilon\cap(V\times\overline X)$
such that $f_J(z,w)$ is holomorphic in $w$ (i.e.\ $\dbar_w f=0$)
and both $f$ and $\dbar f=\dbar_zf$ are in $L^2(p^{-1}(K))$ for all
compact subsets $K\compact V$ (here $\dbar$ operators are of course
taken in the sense of distributions). By construction, we get
a complex of sheaves $(\cF_\varepsilon^\bullet,\dbar\,)$ and
the kernel $\Ker\dbar:\cF_\varepsilon^0\to\cF_\varepsilon^1$ coincides with $\cB_\varepsilon$.
In that sense, if we define $\cO_{L^2}(B_\varepsilon)$ to be the sheaf of
$L^2_\loc$
sections $f$ of $B_\varepsilon$ such that $\dbar f=0$ in the sense
of distributions, then we exactly have $\cO_{L^2}(B_\varepsilon)=\cB_\varepsilon$
as a sheaf. For $z\in V$, the restriction map
$\cB_\varepsilon(V)=\cO_{L^2}(B_\varepsilon)(V)\to B_{\varepsilon,z}$ is an
unbounded closed operator with dense domain, and the kernel is the closure of
$\gm_z\cB_\varepsilon(V)$, which need not be closed. If one insists on getting
continuous fiber restrictions, one could consider the subsheaf
$$
\cO_{C^k}(B_\varepsilon)(V):=
\cO_{L^2}(B_\varepsilon)(V)\cap\cC^k(B_\varepsilon)(V)
$$
where $\cC^k(B_\varepsilon)$ is the sheaf of sections $f$ such that
$\nabla^\ell f$ is continuous in the Hilbert bundle topology for all
real analytic connections $\nabla$ on $B_\varepsilon$ and all
$\ell=0,1,\ldots,k$. For these subsheaves (and any $k\ge 0$), we do get
continuous fiber restrictions $\cO_{C^k}(B_\varepsilon)(V)\to B_{\varepsilon,z}$
for~$z\in V$. In the same way, we could introduce the Dolbeault
complex $\cF_\varepsilon^\bullet\cap\cC^\infty$ and check
that it is a resolution of $\cO\cap \cC^\infty(B_\varepsilon)$,
but we will not need this refinement. Howe\-ver, a useful observation
is that the closed and densely defined operator
$\cO_{L^2}(B_\varepsilon)(V)\to B_{\varepsilon,z}$ is surjective,
in fact it is even true that $H^0(X,\cO_{L^2}(B_\varepsilon))\to
B_{\varepsilon,z}$ is surjective by the Ohsawa-Takegoshi extension
theorem [OhT87] applied on the Stein manifold $U_\varepsilon$. We are
going to see that $\cB_\varepsilon$ can somehow be seen as an
infinite dimensional very ample sheaf. This is already illustrated
by the following result.

\claim 2.5. Proposition|Assume here that $\varepsilon>0$ is taken so small
that $\psi(z,w):=|\logh_z(w)|^2$ is strictly plurisubharmonic up to the
boundary on the compact set $\overline U_\varepsilon
\subset X\times\overline X$. Then the complex of sheaves
$(\cF_\varepsilon^\bullet,\dbar)$ is a resolution of $\cB_\varepsilon$ by soft
sheaves over~$X$ $($actually, by $\cC^\infty_X$-modules$\,)$,
and for every holomorphic vector bundle
$E\to X$ and every $q\geq 1$ we have
$$
H^q(X,\cB_\varepsilon\otimes\cO(E))=H^q\big(\Gamma(X,
\cF_\varepsilon^\bullet\otimes\cO(E)),\dbar\big)=0.
$$
Moreover the fibers $B_{\varepsilon,z}\otimes E_z$ are always generated
by global sections of $H^0(X,\cB_\varepsilon\otimes\cO(E))$, in the sense
that $H^0(X,\cB_\varepsilon\otimes\cO(E))\to B_{\varepsilon,z}\otimes E_z$
is a closed and densely defined operator with surjective image.
\endclaim

\proof{Proof} By construction, we can equip $U_\varepsilon$ with the
the associated K\"ahler metric $\omega=\ii\ddbar\psi$ which is smooth
and strictly positive on $\overline U_\varepsilon$.
We can then take an arbitrary smooth hermitian metric $h_E$ on $E$ and
multiply it by $e^{-C\psi}$, $C\gg 1$, to obtain a bundle with arbitrarily
large positive curvature tensor. The exactness of $\cF_\varepsilon^\bullet$
and cohomology vanishing then follow from
the standard H\"ormander $L^2$ estimates applied either 
locally on $p^{-1}(V)$ for small Stein open sets $V\subset X$,
or globally on $U_\varepsilon$. The global generation of fibers is again
a consequence of the Ohsawa-Takegoshi $L^2$ extension theorem.\qed
\medskip

\claim 2.6. Remark|{\rm
The same result holds for an
arbitrary coherent sheaf $\cE$ instead of a locally free sheaf $\cO(E)$,
the reason being that $p^*\cE$ admits a resolution by (finite dimensional)
locally free sheaves $\cO_{U_{\varepsilon'}}^{\oplus N}$ on a Stein
neighborhood $U_{\varepsilon'}$ of~$\overline U_\varepsilon$.}
\endclaim

\claim 2.7. Remark|{\rm
A strange consequence of these results is that we get some sort
of ``holomorphic embedding''
of an arbitrary complex manifold $X$ into a ``Hilbert Grassmannian'',
mapping every point $z\in X$ to the
closed subspace $S_z$ in the Hilbert space $\bH=\cB_\varepsilon(X)$,
consisting of sections 
$f\in\bH$ such that $f(z)=0$ in $B_{\varepsilon,z}$, i.e.\ $f_{|p^{-1}(z)}=0$.
However, the fact that the restriction morphisms $f\mapsto f_{|p^{-1}(z)}$
are not continuous in $L^2$ norm implies that the map $z\mapsto S_z$ is
not even continuous in the strong topology, i.e.\ the metric topology
for which the distance of two fibers $S_{z_1}$, $S_{z_2}$ is the 
Hausdorff distance of their unit balls in the
$L^2$ norm of~$\cB_\varepsilon(X)$.}
\endclaim

\plainsection{3}{Curvature tensor of Bergman bundles}

\plainsubsection 3.A. Calculation in the model case ($\bC^n$, std)|

In the model situation $X=\bC^n$ with its standard hermitian metric,
we consider the tubular neighborhood
$$
U_\varepsilon:=\{(z,w)\in\bC^n\times\bC^n\,;\;
|\overline w-z|<\varepsilon\}
\leqno(3.1)
$$
and the projections
$$
p=(\pr_1)_{|U_\varepsilon}:U_\varepsilon\to X=\bC^n,~~
(z,w)\mapsto z,\quad
\overline p=(\pr_2)_{|U_\varepsilon}:U_\varepsilon\to X=\bC^n,~~
(z,w)\mapsto w
$$
If one insists on working on a compact complex manifold, the geometry
is locally identical to that of a complex torus $X=\bC^n/\Lambda$ equipped
with a constant hermitian metric $\gamma$. 

\claim 3.2. Remark|{\rm We check here that the Bergman bundle
$B_\varepsilon$ is not holomorphically locally
trivial, even in the above situation where we have invariance by translation.
However, in the category of real analytic bundles, there is a global
trivialization of $B_\varepsilon\to\bC^n$ given by the map
$$
\tau:B_\varepsilon\mathop{\longrightarrow}\limits^{\simeq}
\bC^n\times\cH^2(\bB_n),\quad B_{\varepsilon,z}\ni f_z\longmapsto\tau(f)=(z,g_z),
\quad g_z(w):=f_z(\varepsilon w+\overline z),~~w\in\bB_n,
$$
in other words, for any open set $V\subset\bC^n$ and any $k\in\bN\cup
\{\infty,\omega\}$, we have isomorphisms
$$
C^k(V,B_\varepsilon)\to C^k(V,\cH^2(\bB^n)),\quad
f\mapsto g,\quad g(z,w)=f(z,\varepsilon w+\overline z),
$$
where $f,g$ are $C^k$ in $(z,w)$, holomorphic in $w$, and the derivatives
$z\mapsto D^\alpha_zg(z,\bu)$, $|\alpha|\le k$, define continuous
maps $V\to\cH^2(\bB_n)$.
The complex structures of these bundles are defined by the
$(0,1)$-connections $\dbar_z$ of the associated Dolbeault complexes,
but obviously $\dbar_zf$~and $\dbar_zg$ do not match. In fact,
if we write
$$
g(z,w)=u(z,w)\,dw_1\wedge\ldots\wedge dw_n\in
C^\infty(V,\cH^2(\bB_n))=C^\infty(V)\mathop{\widehat\otimes}\cH^2(\bB_n)
$$
where $\widehat\otimes$ is the $\varepsilon$ or $\pi$-topological
tensor product in the sense of [Gro55], we get
$$
\plaineqalign{
&f(z,w)=g(z,(w-\overline z)/\varepsilon)=
\varepsilon^{-n}u(z,(w{-}\overline z)/\varepsilon)\,
dw_1\wedge\ldots\wedge dw_n,\cr
&\dbar_zf(z,w)=\varepsilon^{-n}\Big(
\dbar_zu(z,(w{-}\overline z)/\varepsilon)-\varepsilon^{-1}\!\!
\sum_{1\leq j\leq n}{\partial u\over\partial w_j}
(z,(w{-}\overline z)/\varepsilon)\,d\overline z_j\Big)\wedge
dw_1\wedge\ldots\wedge dw_n.\cr}
$$
Therefore the trivialization $\tau_*:f\mapsto u$ yields at the level of
$\dbar$-connections an identification
$$
\tau_*:\dbar_zf\mathop{\longmapsto}\limits^\simeq
\dbar_zu+Au
$$
where the ``connection matrix''
$A\in \Gamma(V,\Lambda^{0,1}T^*_X\otimes_\bC
\End(\cH^2(\bB_n)))$
is the constant unbounded Hilbert space operator $A(z)=A$ given by
$$
A:\cH^2(\bB_n)\to\Lambda^{0,1}T^*_X\otimes_\bC\cH^2(\bB_n),\quad
u\mapsto Au=
-\varepsilon^{-1}\sum_{1\leq j\leq n}{\partial u\over \partial w_j}\,
d\overline z_j.
$$
We see that the holomorphic structure of $B_\varepsilon$ is given
by a $(0,1)$-connection that differs by the matrix $A$ from the trivial
$(0,1)$-connection, and as $A$ is unbounded, there is no way we can
make it trivial by a real analytic gauge change with values in Lie
algebra of continuous endomorphisms of $\cH^2(\bB_n)$.\qed
}
\endclaim

We are now going to compute the curvature tensor of the Bergman bundle
$B_\varepsilon$. For the sake of simplicity, we identify here
$\cH^2(\bB_n)$ to the Hardy space of $L^2$ holomorphic functions via
$u\mapsto g=u(w)\,dw_1\wedge\ldots\wedge dw_n$. After rescaling,
we can also assume $\varepsilon=1$, and at least in a first step,
we perform our calculations on $B_1$ rather than~$B_\varepsilon$.
Let us write $w^\alpha=\prod_{1\leq j\leq n}w_j^{\alpha_j}$ for a multiindex
$\alpha=(\alpha_1,\ldots,\alpha_n)\in\bN^n$, and denote by $\lambda$ the
Lebesgue measure on~$\bC^n$. A~well known calculation gives
$$
\int_{\bB_n}|w^\alpha|^2d\lambda(w)=\pi^n{\alpha_1!\ldots\alpha_n!\over(|\alpha|+n)!},
\quad |\alpha|=\alpha_1+\cdots+\alpha_n.
$$
In fact, by using polar coordinates $w_j=r_je^{i\theta_j}$ and
writing $t_j=r_j^2$, we get
$$
\int_{\bB_n}|w^\alpha|^2d\lambda(w)=(2\pi)^n\int_{r_1^2+\cdots+r_n^2<1}
r^{2\alpha}\,r_1dr_1\ldots r_ndr_n=\pi^n I(\alpha)
$$
with
$$
I(\alpha)=\pi^n\int_{t_1+\cdots+t_n<1}
t^\alpha\,dt_1\ldots dt_n.
$$
Now, an induction on $n$ together with the Fubini formula gives
$$
\plaineqalign{
I(\alpha)&=\int_0^1t_n^{\alpha_n}dt_n
\int_{t_1+\cdots+t_{n-1}<1-t_n}(t')^{\alpha'}dt_1\ldots dt_{n-1}\cr
&=I(\alpha')\int_0^1(1-t_n)^{\alpha_1+\cdots+\alpha_{n-1}+n-1}t_n^{\alpha_n}
dt_n\cr}
$$
where $t'=(t_1,\ldots,t_{n-1})$ and $\alpha'=(\alpha_1,\ldots,\alpha_{n-1})$.
As $\int_0^1x^a(1-x)^bdt={a!b!\over(a+b+1)!}$, we get inductively
$$
I(\alpha)={(|\alpha'|+n-1)!\alpha_n!\over(|\alpha|+n)!}I(\alpha')~~~
\Rightarrow~~~I(\alpha)={\alpha_1!\ldots\alpha_n!\over(|\alpha|+n)!}.
$$
Such formulas were already used by Shiffman and Zelditch [ShZ99] in their
study of zeros of random sections of positive line bundles.
They imply that a Hilbert (orthonormal) basis of $\cO\cap L^2(\bB_n)
\simeq\cH^2(\bB_n)$ is
$$
e_\alpha(w)=\pi^{-n/2}
\sqrt{{(|\alpha|+n)!\over\alpha_1!\ldots\alpha_n!}}\,w^\alpha.
\leqno(3.3)
$$
As a consequence, and quite classically, the Bergman kernel of the unit
ball $\bB_n\subset\bC^n$ is
$$
K_n(w)=\sum_{\alpha\in\bN^n}|e_\alpha(w)|^2
=\pi^{-n}\sum_{\alpha\in\bN^n}
{(|\alpha|+n)!\over\alpha_1!\ldots\alpha_n!}\,|w^\alpha|^2
=n!\,\pi^{-n}(1-|w|^2)^{-n-1}.\leqno(3.4)
$$
If we come back to $U_\varepsilon$ for $\varepsilon>0$ not necessarily equal
to $1$ (and do not omit any more the trivial $n$-form
$dw_1\wedge\ldots\wedge dw_n$), we have
to use a rescaling $(z,w)\mapsto(\varepsilon^{-1}z,\varepsilon^{-1}w)$.
This gives for the Hilbert bundle $B_\varepsilon$ a real analytic
orthonormal frame
$$
e_\alpha(z,w)=\pi^{-n/2}\varepsilon^{-|\alpha|-n}
\sqrt{{(|\alpha|+n)!\over\alpha_1!\,...\,\alpha_n!}}\;
(w-\overline z)^\alpha\;dw_1\wedge\ldots\wedge dw_n
\leqno(3.5)
$$
A germ of holomorphic section $\sigma\in\cO_{L^2}(B_\varepsilon)$ near
$z=0$ (say) is thus given by a convergent power series
$$
\sigma(z,w)=\sum_{\alpha\in \bN^n}\xi_\alpha(z)\,e_\alpha(z,w)
$$
such that the functions $\xi_\alpha$ are real analytic on a neighborhood of $0$
and satisfy the following two conditions:
$$\plainleqalignno{
&|\sigma(z)|_h^2:=\sum_{\alpha\in \bN^n}|\xi_\alpha(z)|^2~~\hbox{converges in $L^2$ near $0$},
&(3.6)\cr
&\dbar_{z_k}\sigma(z,w)=\sum_{\alpha\in \bN^n}
\dbar_{z_k}\xi_\alpha(z)\,e_\alpha(z,w)+
\xi_\alpha(z)\,\dbar_{z_k}e_\alpha(z,w)\equiv 0.&(3.7)\cr}
$$
Let $c_k=(0,\ldots,1,\ldots, 0)$ be the canonical basis of
the $\bZ$-module $\bZ^n$. A straightforward calculation from (3.5) yields
$$
\dbar_{z_k}e_\alpha(z,w)=-\varepsilon^{-1}
\sqrt{\alpha_k(|\alpha|+n)}\; e_{\alpha-c_k}(z,w).
$$
We have the slight problem that the coefficients are unbounded as
$|\alpha|\to+\infty$, and therefore the two terms occurring in
(3.7) need not form convergent series when taken separately. However
if we take $\sigma\in\cO_{L^2}(B_{\varepsilon'})$ in a slightly bigger tubular
neighborhood ($\varepsilon'>\varepsilon$), the $L^2$ condition implies
that $\sum_\alpha(\varepsilon''/\varepsilon)^{2|\alpha|}|\xi_\alpha|^2$ is
uniformly convergent for every $\varepsilon''\in{}]\varepsilon,\varepsilon'[\,$,
and this is more than enough to ensure convergence, since the growth of
$\alpha\mapsto\sqrt{\alpha_k(|\alpha|+n)}$ is at most linear; we can even
iterate as many derivatives as we want. For a smooth section
$\sigma\in C^\infty(B_{\varepsilon'})$, the coefficients
$\xi_\alpha$ are smooth, with $\sum(\varepsilon'/\varepsilon)^{2|\alpha|}
|\partial_z^\beta\dbar_z^\gamma\xi_\alpha|^2$
convergent for all $\beta,\gamma$, and we get
$$
\plaineqalign{
\dbar_{z_k}\sigma(z,w)
&=\sum_{\alpha\in\bN^n}\dbar_{z_k}\xi_\alpha(z)\,e_\alpha(z,w)+
\xi_\alpha(z)\,\dbar_{z_k}e_\alpha(z,w)\cr
&=\sum_{\alpha\in\bN^n}\dbar_{z_k}\xi_\alpha(z)\,e_\alpha(z,w)-
\varepsilon^{-1}\sqrt{\alpha_k(|\alpha|+n)}\;
\xi_\alpha(z)\;e_{\alpha-c_k}(z,w)\cr
&=\sum_{\alpha\in\bN^n}\big(\dbar_{z_k}\xi_\alpha(z)
-\varepsilon^{-1}\sqrt{(\alpha_k+1)(|\alpha|+n+1)}\;
\xi_{\alpha+c_k}(z)\big)\;e_\alpha(z,w),\cr}
$$
after replacing $\alpha$ by $\alpha+c_k$ in the terms containing
$\varepsilon^{-1}$. The $(0,1)$-part $\nabla_h^{0,1}$ of the Chern connection
$\nabla_h$ of $(B_\varepsilon,h)$ with respect to the orthonormal
frame~$(e_\alpha)$ is thus given by
$$
\nabla_h^{0,1}\sigma
=\sum_{\alpha\in\bN^n}\Big(\dbar\xi_\alpha
-\sum_k\varepsilon^{-1}\sqrt{(\alpha_k+1)(|\alpha|+n+1)}\;
\xi_{\alpha+c_k}\;d\overline z_k\Big)\otimes e_\alpha.\leqno(3.8)
$$
The $(1,0)$-part can be derived from the identity
$\partial|\sigma|_h^2=\langle\nabla_h^{1,0}\sigma,\sigma\rangle_h+
\langle\sigma,\nabla_h^{0,1}\sigma\rangle_h$.
However
$$
\plaineqalign{
\partial_{z_j}|\sigma|_h^2
&=\partial_{z_j}\sum_{\alpha\in\bN^n}\xi_\alpha\overline\xi_\alpha=
\sum_{\alpha\in\bN^n}(\partial_{z_j}\xi_\alpha)\;\overline\xi_\alpha+
\xi_\alpha\;(\overline{\dbar_{z_j}\xi_\alpha}\,)\cr
&=\sum_{\alpha\in\bN^n}
\Big(~\partial_{z_j}\xi_\alpha+\varepsilon^{-1}
\sqrt{\alpha_j(|\alpha|+n)}\;
\xi_{\alpha-c_j}\Big)\,\overline\xi_\alpha\cr
&+\sum_{\alpha\in\bN^n}\xi_\alpha\Big(~\overline{\dbar_{z_j}\xi_\alpha
-\varepsilon^{-1}\sqrt{(\alpha_j+1)(|\alpha|+n+1)}\;
\xi_{\alpha+c_j}}~\Big).\cr}
$$
For $\sigma\in C^\infty(B_{\varepsilon'})$, it follows from there that 
$$
\nabla_h^{1,0}\sigma=
\sum_{\alpha\in\bN^n}\Big(\partial\xi_\alpha+\varepsilon^{-1}\sum_j
\sqrt{\alpha_j(|\alpha|+n)}\;\xi_{\alpha-c_j}dz_j\Big)\otimes e_\alpha.
\leqno(3.9)
$$
Finally, to find the curvature tensor of $(B_\varepsilon,h)$, we only
have to compute the $(1,1)$-form $(\nabla_h^{1,0}\nabla_h^{0,1}+
\nabla_h^{0,1}\nabla_h^{1,0})\sigma$ and take the terms that contain
no differentiation at all, especially in view of the usual identity
$\partial\dbar+\dbar\partial=0$ and the fact that we also have here
$(\nabla_h^{1,0})^2=0$, $(\nabla_h^{0,1})^2=0$. As
$(\alpha-c_j)_k=\alpha_k-\delta_{jk}$ and $(\alpha+c_k)_j=
\alpha_j+\delta_{jk}$, we are left with
$$
\plaineqalign{
&\big(\nabla_h^{1,0}\nabla_h^{0,1}+\nabla_h^{0,1}\nabla_h^{1,0}\big)\sigma\cr
&~~{}=
-\,\varepsilon^{-2}\sum_{\alpha\in\bN^n}\sum_{j,k}
\sqrt{\alpha_j(|\alpha|{+}n)}\;
\sqrt{(\alpha_k{-}\delta_{jk}{+}1)(|\alpha|{+}n)}\;
\xi_{\alpha-c_j+c_k}\;dz_j\wedge d\overline z_k\otimes e_\alpha\cr
&~~{\phantom{{}={}}}{+}\,\varepsilon^{-2}\sum_{\alpha\in\bN^n}\sum_{j,k}
\sqrt{(\alpha_j{+}\delta_{jk})(|\alpha|{+}n{+}1)}\;
\sqrt{(\alpha_k{+}1\rlap{\hbox{$\phantom{\alpha_j}$}})(|\alpha|{+}n{+}1)}\;
\xi_{\alpha-c_j+c_k}\;dz_j\wedge d\overline z_k\otimes e_\alpha\cr
&~~{}=
-\,\varepsilon^{-2}\sum_{\alpha\in\bN^n}\sum_{j,k}
\sqrt{(\alpha_j{-}\delta_{jk})(\alpha_k{-}\delta_{jk})}\;(|\alpha|+n-1)\;
\xi_{\alpha-c_j}\;dz_j\wedge d\overline z_k\otimes e_{\alpha-c_k}\cr
&~~{\phantom{{}={}}}{+}\,\varepsilon^{-2}\sum_{\alpha\in\bN^n}\sum_{j,k}
\sqrt{\alpha_j\alpha_k}\;(|\alpha|+n)\;
\xi_{\alpha-c_j}\;dz_j\wedge d\overline z_k\otimes e_{\alpha-c_k}.\cr
&~~{}=\phantom{{}-\,}
\varepsilon^{-2}\sum_{\alpha\in\bN^n}\sum_{j,k}
\sqrt{\alpha_j\alpha_k}\;\xi_{\alpha-c_j}\;
dz_j\wedge d\overline z_k\otimes e_{\alpha-c_k}\cr
&~~{\phantom{{}={}}}{+}\,\varepsilon^{-2}\sum_{\alpha\in\bN^n}
\sum_{j}(|\alpha|+n-1)\;
\xi_{\alpha-c_j}\;dz_j\wedge d\overline z_j\otimes e_{\alpha-c_j},\cr}
$$
where the last summation comes from the subtraction of the diagonal terms
$j=k$. By changing $\alpha$ into $\alpha+c_j$ in that summation,
we obtain the following expression of the curvature tensor
of $(B_\varepsilon,h)$.

\claim 3.10. Theorem|The curvature tensor of the Bergman bundle
$(B_\varepsilon,h)$ is given by
$$
\langle\Theta_{B_\varepsilon,h}\sigma(v,Jv),\sigma\rangle_h
=\varepsilon^{-2}\sum_{\alpha\in\bN^n}\Bigg(\bigg|
\sum_j\sqrt{\alpha_j}\;\xi_{\alpha-c_j}v_j\bigg|^2
+\sum_{j}(|\alpha|+n)\;|\xi_\alpha|^2|v_j|^2\Bigg)
$$
for every $\sigma=\sum_\alpha\xi_\alpha e_\alpha\in B_{\varepsilon'}$,
$\varepsilon'>\varepsilon$, and every
tangent vector $v=\sum v_j\;\partial/\partial z_j$.
\endclaim

The above curvature hermitian tensor is positive definite, and even
positive definite unbounded if we view it as a hermitian form
on $T_X\otimes B_\varepsilon$ rather than on $T_X\otimes B_{\varepsilon'}$.
This is not so surprising since the connection matrix was already an
unbounded operator. Philosophically, the very ampleness of the
sheaf~$\cB_\varepsilon$
was also a strong indication that the curvature of the corresponding
vector bundle $B_\varepsilon$ should have been positive.
Observe that we have in fact
$$
\plainleqalignno{
\varepsilon^{-2}\sum_{\alpha\in\bN^n}
\sum_{j}&(|\alpha|+n)\;|\xi_\alpha|^2|v_j|^2\cr
&\leq\langle\Theta_{B_\varepsilon,h}\sigma(v,Jv),\sigma\rangle_h
\leq 2\varepsilon^{-2}\sum_{\alpha\in\bN^n}
\sum_{j}(|\alpha|+n)\;|\xi_\alpha|^2|v_j|^2,&(3.11)\cr}
$$
thanks to the Cauchy-Schwarz inequality
$$
\plaineqalign{
\sum_{\alpha\in\bN^n}\bigg|\sum_j\sqrt{\alpha_j}\;\xi_{\alpha-c_j}v_j\bigg|^2
&\leq\sum_\ell|v_\ell|^2\sum_{\alpha\in\bN^n}\sum_j\alpha_j|\xi_{\alpha-c_j}|^2
=\sum_\ell|v_\ell|^2\sum_j\sum_{\alpha\in\bN^n}\alpha_j|\xi_{\alpha-c_j}|^2\cr
&=\sum_\ell|v_\ell|^2\sum_j\sum_{\alpha\in\bN^n}(\alpha_j+1)|\xi_\alpha|^2
=\sum_\ell|v_\ell|^2\sum_{\alpha\in\bN^n}(|\alpha|+n)|\xi_\alpha|^2.\cr}
$$

\plainsubsection 3.B. Curvature of Bergman bundles on compact
hermitian manifolds|

We consider here the general situation of a compact hermitian manifold
$(X,\gamma)$ described in \S1, where $\gamma$ is real analytic and
$\exph$ is the associated partially holomorphic exponential map.
Fix a point $x_0\in X$, and use a holomorphic system of coordinates
$(z_1,\ldots,z_n)$ centered at $x_0$, provided by
$$
\exph_{x_0}:T_{X,x_0}\supset V\to X.
$$
If we take $\gamma_{x_0}$ orthonormal coordinates on $T_{X,x_0}$, then by
construction the fiber of $p:U_\varepsilon\to X$ over $x_0$ is the
standard $\varepsilon$-ball in the coordinates $(w_j)=(\overline z_j)$.
Let $T_X\to V\times\bC^n$ be the trivialization of $T_X$ in the
coordinates~$(z_j)$, and
$$
X\times X\to T_X,\quad
(z,w)\mapsto\xi=\logh_z(w)
$$
the expression of $\logh$ near $(x_0,\overline x_0)$,
that is, near $(z,w)=(0,0)$. By our choice of coordinates, we have
$\logh_0(w)=w$ and of course $\logh_z(z)=0$, hence
we get a real analytic expansion of the form
$$
\plaineqalign{
\logh_z(w)=w-z&+\sum z_ja_j(w\,{-}\,z)+\sum\overline z_ja'_j(w\,{-}\,z)\cr
&+\sum z_jz_kb_{jk}(w\,{-}\,z)
+\sum \overline z_j\overline z_kb'_{jk}(w\,{-}\,z)
+\sum z_j\overline z_kc_{jk}(w\,{-}\,z)+O(|z|^3)\cr}
$$
with holomorphic coefficients $a_j$, $a'_j$, $b_{jk}$, $b'_{jk}$, $c_{jk}$
vanishing at~$0$. In fact by [Dem94], we always have
$da'_j(0)=0$, and if $\gamma$ is K\"ahler, the equality $da_j(0)=0$
also holds; we will not use these properties  here. In coordinates, we
then have locally near $(0,0)\in\bC^n\times\bC^n$
$$
U_{\varepsilon,z}=\big\{(z,w)\in\bC^n\times\bC^n\,;\;
|\Psi_z(w)|<\varepsilon\big\}
$$
where $\Psi_z(w)=\overline{\logh_z(\overline w)}$ has a similar expansion
$$
\plainleqalignno{
\Psi_z(w)=w-\overline z&+\sum z_ja_j(w\,{-}\,\overline z)
+\sum\overline z_ja'_j(w\,{-}\,\overline z)\cr
&+\sum z_jz_kb_{jk}(w\,{-}\,\overline z)
+\sum \overline z_j\overline z_kb'_{jk}(w\,{-}\,\overline z)
+\sum z_j\overline z_kc_{jk}(w\,{-}\,\overline z)+O(|z|^3)&(3.12)\cr}
$$
(when going from $\logh$ to $\Psi$, the coefficients $a_j$, $a'_j$ and
$b_j$, $b'_j$ get twisted, but we do not care and keep the same
notation for $\Psi$, as we will not refer to $\logh$ any more).
In this situation, the Hilbert bundle $B_\varepsilon$
has a real analytic normal
frame given by $\wt e_\alpha=\Psi^*e_\alpha$ where
$$
e_\alpha(w)=\pi^{-n/2}\varepsilon^{-|\alpha|-n}
\sqrt{{(|\alpha|+n)!\over\alpha_1!\,...\,\alpha_n!}}\;
w^\alpha\;dw_1\wedge\ldots\wedge dw_n
\leqno(3.13)
$$
and the pull-back $\Psi^*e_\alpha$ is taken with respect to $w\mapsto\Psi_z(w)$
($z$ being considered as a parameter). For a local section
$\sigma=\sum_\alpha\xi_\alpha\wt e_\alpha\in C^\infty(B_{\varepsilon'})$,
$\varepsilon>\varepsilon$, we can write
$$
\dbar_{z_k}\sigma(z,w)=
\sum_{\alpha\in\bN^n}\dbar_{z_k}\xi_\alpha(z)\,\wt e_\alpha(z,w)+
\xi_\alpha(z)\,\dbar_{z_k}\wt e_\alpha(z,w).
$$
Near $z=0$, by taking the derivative of $\Psi^*e_\alpha(z,w)$, we find
$$
\plaineqalign{
\dbar_{z_k}\wt e_\alpha(z,w)=
&-\varepsilon^{-1}\sqrt{\alpha_k(|\alpha|+n)}\;\wt e_{\alpha-c_k}(z,w)\cr
&+\varepsilon^{-1}\sum_m
\sqrt{\alpha_m(|\alpha|+n)}\;
\bigg(a'_{k,m}(w\,{-}\,\overline z)+
\sum_jz_jc_{jk,m}(w\,{-}\,\overline z)\bigg)\;\wt e_{\alpha-c_m}(z,w)\cr
&+\sum_m \bigg(
{\partial a'_{k,m}\over\partial w_m}(w\,{-}\,\overline z)+
\sum_jz_j{\partial c_{jk,m}\over\partial w_m}(w\,{-}\,\overline z)\bigg)
\;\wt e_\alpha(z,w)
+O(\overline z,|z|^2),\cr}
$$
where the last  sum comes from the expansion of
$dw_1\wedge\ldots\wedge dw_n$, and $a'_{k,m}$, $c_{jk,m}$ are
the $m$-th components of~$a'_k$ and $c_{jk}$. This
gives two additional terms in comparison to the translation invariant case,
but these terms are ``small'' in the sense that the first one vanishes
at $(z,w)=(0,0)$ and the second one does not involve~$\varepsilon^{-1}$.
If $\nabla_{h,0}^{0,1}$ is the \hbox{$\dbar$-connection} associated with the
standard tubular neighborhood $|\overline w-z|<\varepsilon$, we thus find
in terms of the local trivialization $\sigma\simeq\xi=\sum\xi_\alpha
\wt\varepsilon_\alpha$ an expression of the form
$$
\nabla_h^{0,1}\sigma\simeq\nabla_{h,0}^{0,1}\xi+A^{0,1}\xi,
$$
where
$$
\plaineqalign{
A^{0,1}\bigg(\sum_\alpha&\xi_\alpha\wt e_\alpha\bigg)
=\sum_{\alpha\in\bN^n}\sum_k\cr
&\xi_\alpha\Bigg(\varepsilon^{-1}\sum_m
\sqrt{\alpha_m(|\alpha|+n)}\;
\bigg(a'_{k,m}(w)+\sum_jz_jc_{jk,m}(w)\bigg)\;
d\overline z_k\otimes\wt e_{\alpha-c_m}\cr
&\qquad{}+\sum_m \bigg(
{\partial a'_{k,m}\over\partial w_m}(w)+
\sum_jz_j{\partial c_{jk,m}\over\partial w_m}(w)\bigg)
\;d\overline z_k\otimes\wt e_\alpha\Bigg)+O(\overline z,|z|^2).\cr}
$$
The corresponding $(1,0)$-parts satisfy
$$
\nabla_h^{1,0}\sigma\simeq\nabla_{h,0}^{1,0}\xi+A^{1,0}\xi,\qquad
A^{1,0}=-(A^{0,1})^*,
$$
and the corresponding curvature tensors are related by
$$
\Theta_{\beta_\varepsilon,h}=
\Theta_{\beta_\varepsilon,h,0}+\partial A^{0,1}+\dbar A^{1,0}+
A^{1,0}\wedge A^{0,1}+A^{0,1}\wedge A^{1,0}.\leqno(3.14)
$$
At $z=0$ we have
$$
\plaineqalign{
A^{0,1}\xi
=\sum_{\alpha\in\bN^n}\sum_k\xi_\alpha\Bigg(\varepsilon^{-1}
&\sum_{m}
\sqrt{\alpha_m(|\alpha|+n)}\;a'_{k,m}(w)\;
d\overline z_k\otimes\wt e_{\alpha-c_m}\cr
+&\sum_m{\partial a'_{k,m}\over\partial w_m}(w)\;
d\overline z_k\otimes\wt e_\alpha\Bigg),\cr
\partial A^{0,1}\xi
=\sum_{\alpha\in\bN^n}\sum_k\xi_\alpha\Bigg(\varepsilon^{-1}
&\sum_{j,m}
\sqrt{\alpha_m(|\alpha|+n)}\;c_{jk,m}(w)\;dz_j\wedge
d\overline z_k\otimes\wt e_{\alpha-c_m}\cr
+&\sum_{j,m}{\partial c_{jk,m}\over\partial w_m}(w)\;
dz_j\wedge d\overline z_k\otimes\wt e_\alpha\Bigg),\cr}
$$
and $A^{1,0}$, $\dbar A^{1,0}$ are, up to the sign, the adjoint endomorphisms
of $A^{0,1}$ and $\partial A^{0,1}$. The unboundedness comes from the fact
that we have unbounded factors $\sqrt{(\alpha_m+1)(|\alpha|+n+1)}\,$;
it is worth noticing that multiplication by a holomorphic factor $u(w)$ is
a continuous operator on the fibers $B_{\varepsilon,z}$, whose norm remains
bounded as $\varepsilon\to 0$. In this setting, it can be seen that
the only term in (3.14) that is (a priori) not small with respect to the
main term $\Theta_{B_\varepsilon,h,0}$ is the term involving
$\varepsilon^{-2}$ in $A^{1,0}\wedge A^{0,1}+A^{0,1}\wedge A^{1,0}$,
and that the other terms appearing in the quadratic
form~$\langle\Theta_{B_\varepsilon,h}\xi,\xi\rangle$
are $O(\varepsilon^{-1}\sum(|\alpha|+n)|\xi_\alpha|^2)$ or smaller.
In order to check this, we expand $c_{jk,m}(w)$ into a power
series $\sum_\mu c_{jk,m,\mu}\;g_\mu(w)$
where
$$
g_\mu(w)=s_\mu^{-1}w^\mu,\quad\hbox{with
$\displaystyle
s_\mu=\sup_{|w|\leq 1}|w^\mu|
=\prod_{1\leq j\leq n}\bigg({\mu_j\over|\mu|}\bigg)^{\mu_j/2}
={\prod\mu_j^{\mu_j/2}\over|\mu|^{|\mu|/2}}$},
\leqno(3.15)
$$
so that $\sup_{|w|\leq\varepsilon}|g_\mu(w)|=\varepsilon^{|\mu|}$.
We get from the term
$\langle\partial A^{0,1}\xi,\xi\rangle$ a summation
$$
\Sigma(\xi)=
\varepsilon^{-1}\sum_{j,k,m}~\sum_{\alpha\in\bN^n}\sqrt{\alpha_m(|\alpha|+n)}\;
\sum_{\mu\in\bN^n}c_{jk,m,\mu}\;dz_j\wedge d\overline z_k\otimes
\langle \xi_\alpha\;g_\mu\wt e_\alpha,\xi\rangle.
$$
At $z=0$, $g_\mu\wt e_\alpha=g_\mu e_\alpha$ is proportional to
$e_{\alpha+\mu}$, and by (3.15) and the definition of the
$L^2$ norm, we have $\Vert g_\mu\wt e_\alpha\Vert\leq \varepsilon^{|\mu|}$
and $|\langle \xi_\alpha\;g_\mu\wt e_\alpha,\xi\rangle|\leq
\varepsilon^{|\mu|}\;|\xi_\alpha||\xi_{\alpha+\mu}|$. We infer
$$
\big|\Sigma(\xi)\big|\leq
\varepsilon^{-1}\sum_{j,k,m}~\sum_{\alpha\in\bN^n}\sqrt{\alpha_m(|\alpha|+n)}\;
\sum_{\mu\in\bN^n}|c_{jk,m,\mu}|\;\varepsilon^{|\mu|}\;
|\xi_\alpha||\xi_{\alpha+\mu}|.
$$
Let $r$ be the infimum of the radius of convergence of $w\mapsto\Psi_z(w)$
over all~$z\in X$. Then for $\varepsilon<r$ and $r'\in{}]\varepsilon,r[$,
we have a uniform bound $|c_{jk,m,\mu}|\leq C(1/r')^{|\mu|}$, hence
$$
\big|\Sigma(\xi)\big|\leq
C'\varepsilon^{-1}\sum_{\alpha\in\bN^n}\sum_{\mu\in\bN^n}
\Big({\varepsilon\over r'}\Big)^{|\mu|}
\sqrt{\alpha_m(|\alpha|+n)}\;|\xi_\alpha||\xi_{\alpha+\mu}|.
$$
If we write
$$
\plaineqalign{
\sqrt{\alpha_m(|\alpha|+n)}\;|\xi_\alpha||\xi_{\alpha+\mu}|
&\leq{1\over 2}(|\alpha|+n)\big(|\xi_\alpha|^2+|\xi_{\alpha+\mu}|^2\big)\cr
&\leq {1\over 2}\big((|\alpha|+n)|\xi_\alpha|^2+(|\alpha+\mu|+n)
|\xi_{\alpha+\mu}|^2\big),\cr}
$$
the above bound implies
$$
\big|\Sigma(\xi)\big|\leq
C'\varepsilon^{-1}(1-\varepsilon/r')^{-n}
\sum_{\alpha\in\bN^n}(|\alpha|+n)|\xi_\alpha|^2
=O\bigg(\varepsilon^{-1}\sum_{\alpha\in\bN^n}(|\alpha|+n)|\xi_\alpha|^2\bigg).
$$
We now come to the more annoying term
$A^{1,0}\wedge A^{0,1}+A^{0,1}\wedge A^{1,0}$, and especially to the
part containing $\varepsilon^{-2}$ (the other parts can be treated as
above or are smaller). We compute explicitly that term by expanding
$a'_{k,m}(w)$ into a power series $\sum_\mu a'_{k,m,\mu}\;g_\mu(w)$
as above. Let us write $g_m(w)=s_\mu^{-1}w^\mu$.
As $a'_{k,m}(0)=0$, the relevant term in $A^{0,1}$ is
$$
\varepsilon^{-1}\sum_{k,m}\;
\sum_{\mu\in\bN^n\ssm\{0\}}a'_{k,m,\mu}s_\mu^{-1}\;d\overline z_k\otimes W^\mu D_m
$$
where $D_m$ and
$W^\mu=W_1^{\mu_1}\ldots W_n^{\mu_n}$ are
operators on the Hilbert space $\cH^2(B_{\varepsilon,0})$, defined by
$$
D_m\wt e_\alpha=\sqrt{\alpha_m(|\alpha|+n)}\;\wt e_{\alpha-c_m},\quad
W_m(f)=w_m f.
$$
The corresponding term in $A^{1,0}$ is the opposite of the adjoint, namely
$$
-\varepsilon^{-1}\sum_{j,\ell}
\sum_{\lambda\in\bN^n\ssm\{0\}}a'_{k,\ell,\lambda}s_\lambda^{-1}
dz_j\otimes D^*_\ell W^{*\lambda}
$$
and the annoying term in $A^{1,0}\wedge A^{0,1}+A^{0,1}\wedge A^{1,0}$ is
$$
\plainleqalignno{
Q=&-\varepsilon^{-2}\sum_{j,k,\ell,m}
\sum_{\lambda,\mu\in\bN^n\ssm\{0\}}a'_{k,\ell,\lambda}s_\lambda^{-1}\;
a'_{k,m,\mu}s_\mu^{-1}\;dz_j\wedge d\overline z_k\otimes{}&(3.16)\cr
&\kern120pt{}
\Big(D^*_\ell W^{*\lambda}W^\mu D_m-W^\mu D_mD^*_\ell W^{*\lambda}\Big).\cr}
$$
We have here
$\Vert W^\mu\Vert\leq s_\mu\;\varepsilon^{|\mu|}$ (as $W^\mu$ is the 
multiplication by $w^\mu=s_\mu\;g_\mu(w)$, and 
\hbox{$|g_\mu|\leq\varepsilon^{|\mu|}$} on $B_{\varepsilon,0}$). The operators
$D^*_\ell$ and $D_m$ are unbounded, but the important point
is that their commutators have substantially better continuity
than what could be expected a priori. We have for instance
$$
\plaineqalign{
&D_m\wt e_\alpha=\sqrt{\alpha_m(|\alpha|+n)}\;\wt e_{\alpha-c_m},\quad
D_\ell^*(\wt e_\alpha)=\sqrt{(\alpha_\ell+1)(|\alpha|+n+1)}\;
\wt e_{\alpha+c_\ell},\cr
&[D^*_\ell,D_m](\wt e_\alpha)=\Big(
\sqrt{(\alpha_\ell+1-\delta_{\ell m})\alpha_m}\;(|\alpha|+n)\cr
&\kern100pt{}
-\sqrt{(\alpha_\ell+1)(\alpha_m+\delta_{\ell m})}\;(|\alpha|+n+1)\,\Big)\;
\wt e_{\alpha+c_\ell-c_m}\cr}
$$
and the coefficient between braces is controlled by $2(|\alpha|+n)$, as
one sees by considering\break separately the two cases $\ell\neq m$, where we get
$-\sqrt{(\alpha_\ell+1)\alpha_m}$, and $\ell=m$, where we get\break
\hbox{$\alpha_\ell(|\alpha|+n)-(\alpha_\ell+1)(|\alpha|+n+1)$}. Therefore
$\Vert [D^*_\ell,D_m](\wt e_\alpha)\Vert\leq 2(|\alpha|+n)$.
We obtain similarly
$$
\plaineqalign{
&W_m(\wt e_\alpha)=\varepsilon\sqrt{{\alpha_m+1\over|\alpha|+n+1}}\;\wt e_{\alpha+c_m},\quad
W_\ell^*(\wt e_\alpha)=\varepsilon\sqrt{{\alpha_\ell\over|\alpha|+n}}\;
\wt e_{\alpha-c_\ell},\cr
&[W^*_\ell,W_m](\wt e_\alpha)=\varepsilon^2\bigg(
{\sqrt{(\alpha_\ell+\delta_{\ell m})(\alpha_m+1)}\over|\alpha|+n+1}-
{\sqrt{\alpha_\ell(\alpha_m+1-\delta_{\ell m})}\over|\alpha|+n}\,\bigg)\;\wt e_{\alpha-c_\ell+c_m},\cr}
$$
and it is easy to see that the coefficient between large braces is bounded
for $\ell\neq m$ by\break
$\sqrt{\alpha_\ell(\alpha_m+1)}/((|\alpha|+n)(|\alpha|+n+1))
\leq(|\alpha|+n)^{-1}$, and for $\ell=m$ we have as well
$$
\bigg|{(\alpha_\ell+1)(|\alpha|+n)-\alpha_\ell(|\alpha|+n+1)\over
(|\alpha|+n)(|\alpha|+n+1)}\bigg|\leq(|\alpha|+n)^{-1}.
$$
Therefore $\Vert [W^*_\ell,W_m](\wt e_\alpha)\Vert\leq
\varepsilon^2(|\alpha|+n)^{-1}$. Finally
$$
\plaineqalign{
&[W^*_\ell,D_m](\wt e_\alpha)=\varepsilon\Bigg(
\sqrt{{(\alpha_\ell-\delta_{\ell m})\alpha_m(|\alpha|+n)\over
|\alpha|+n-1}}
-\sqrt{{\alpha_\ell(\alpha_m-\delta_{\ell m})(|\alpha|+n-1)\over
|\alpha|+n}}\;\Bigg)\;
\wt e_{\alpha-c_\ell-c_m}\cr}
$$
with a coefficient between braces less than $1$, thus
$\Vert [W^*_\ell,D_m](\wt e_\alpha)\Vert\leq\varepsilon$. By
adjunction, the same is true for $[D^*_\ell,W_m]$, and we can summarize
our estimates as follows:
$$
\plaincases{
\Vert W^\mu\Vert\leq s_\mu\varepsilon^{|\mu|},~~
\Vert W^{*\lambda}\Vert\leq s_\lambda\varepsilon^{|\lambda|},~~
\Vert D^*_\ell(\wt e_\alpha)\Vert\leq|\alpha|+n+1,~~
\Vert D_m(\wt e_\alpha)\Vert\leq |\alpha|+n,\cr
\noalign{\vskip6pt}
\Vert [D^*_\ell,D_m](\wt e_\alpha)\Vert\leq 2(|\alpha|+n),\quad
\Vert [W^*_\ell,W_m](\wt e_\alpha)\Vert\leq (|\alpha|+n)^{-1},\cr
\noalign{\vskip6pt}
\Vert [W^*_\ell,D_m](\wt e_\alpha)\Vert\leq\varepsilon,\quad
\Vert [D^*_\ell,W_m](\wt e_\alpha)\Vert\leq\varepsilon\cr}
\leqno(3.17)
$$
Now, we observe that both $D^*_\ell W^{*\lambda}W^\mu D_m(\wt e_\alpha)$ and
$W^\mu D_mD^*_\ell W^{*\lambda}(\wt e_\alpha)$ are multiples
of\break $\wt e_{\alpha+c_\ell-c_m-\lambda+\mu}$. By considering
the second product $W^\mu D_mD^*_\ell W^{*\lambda}$ and
permuting successively its factors $D_mD^*_\ell$, $D_mW^{*\lambda}$,
$W^\mu D^*_\ell$, $W^\mu W^{*\lambda}$, the diffe\-rence with
$D^*_\ell W^{*\lambda}W^\mu D_m$ is expressed
as a sum of $1+|\lambda|+|\mu|+|\lambda||\mu|$ terms involving commutators.
We derive from our estimates (3.17) precise bounds for the image of 
$\wt e_\alpha$ by the commutators.
For instance, when we arrive at $D^*_\ell W^\mu W^{*\lambda}D_m$
and permute $W^\mu W^{*\lambda}$, we go through intermediate steps
$D_\ell W^{*\lambda'}W^{\mu'}W_k W^*_j W^{\mu''}W^{*\lambda''}D_m$
with $\lambda=\lambda'+\lambda''+c_j$, $\mu=\mu'+\mu''+c_k$,
$|\lambda|=|\lambda'|+|\lambda''|+1$, $|\mu|=|\mu'|+|\mu''|+1$,
and have to evaluate the commutators
$$
D_\ell W^{*\lambda'}W^{\mu'}[W^*_j,W_k]W^{\mu''}W^{*\lambda''}D_m(\wt e_\alpha).
$$
By (3.17), the norm of these $|\lambda||\mu|$ terms is bounded by
$$
\plainleqalignno{
&\big((|\alpha|{-}|\lambda|{+}|\mu|{-}1)_+{+}n{+}1\big)
s_{\lambda'}\varepsilon^{|\lambda'|}s_{\mu'}\varepsilon^{|\mu'|}
\big((|\alpha|{-}|\lambda''|{+}|\mu''|)_{+}{+}n\big)^{-1}
s_{\mu''}\varepsilon^{|\mu''|}
s_{\lambda''}\varepsilon^{|\lambda''|}(|\alpha|{+}n)\cr
&\kern60pt{}\leq
s_{\lambda'}s_{\lambda''}s_{\mu'}s_{\mu''}\;\varepsilon^{|\lambda|+|\mu|}\;
{((|\alpha|-|\lambda|+|\mu|-1)_++n+1)(|\alpha|+n)\over(|\alpha|-|\lambda|)_++n}.
&(3.18)\cr}
$$
The remaining commutators are easier, they lead to bounds
$$
\plaincases{
s_\lambda s_\mu\;\varepsilon^{|\lambda|+|\mu|}\;
2((|\alpha|-|\lambda|)_++n)\hfill&\hbox{(once)},\hfill\cr
\noalign{\vskip6pt}
s_{\lambda'}s_{\lambda''}s_\mu\;\varepsilon^{|\lambda|+|\mu|}\;
(|\alpha|-|\lambda|-1)_++n+1\hfill&\hbox{($|\lambda|$ times)},\hfill\cr
\noalign{\vskip6pt}
s_\lambda s_{\mu'}s_{\mu''}\;\varepsilon^{|\lambda|+|\mu|}\;
(|\alpha|+n)\hfill&\hbox{($|\mu|$ times)}.\hfill\cr}\leqno(3.19)
$$
In the final estimates, we will have to bound some combinatorial
factors of the form
$$
{s_{\lambda'}s_{\lambda''}\over s_\lambda}\;
{s_{\mu'}s_{\mu''}\over s_\mu} \quad\hbox{(worst case)},\leqno(3.20)
$$
and we want the ratios $s_{\lambda'}s_{\lambda''}/s_\lambda$
to be as small as possible (clearly they are at least equal to~$1$).
For this, we try to keep the proportions $\lambda'_j/|\lambda'|$,
$\lambda''_j/|\lambda''|$ as close as possible to $\lambda_j/|\lambda|$
by selecting carefully which factor $W^*_\ell$ (and $W_m$) we exchange
at each step. After a permutation of the coordinates, we may
assume than $\lambda_n\geq\max_{j<n}\lambda_j$, hence
$\lambda_n\geq{1\over n}|\lambda|$. 
If $t'=|\lambda'|/|\lambda|$ and $t''=|\lambda''|/|\lambda|
=1-t'-1/|\lambda|$, we take
$\lambda'_j=\lfloor t'\lambda_j\rfloor$,
$\lambda''_j=\lfloor t''\lambda_j\rfloor$
for $j\leq n-1$ and compensate by taking ad hoc values of $\lambda'_n$,
$\lambda''_n$ and $c_\ell=\lambda-(\lambda'+\lambda'')$. Then
$t'\lambda_j-1<\lambda'_j\leq t'\lambda_j$ for $j<n$ and
$$
\lambda'_n=|\lambda'|-\sum_{j<n}\lambda'_j~~
\plaincases{
\kern-3pt{}<t'|\lambda|-\sum_{j<n}t'\lambda_j+n-1=t'\lambda_n+n-1,\cr
\noalign{\vskip7pt}
\kern-3pt{}\geq t'|\lambda|-\sum_{j<n}t'\lambda_j=t'\lambda_n.\cr}
$$
Therefore 
$$
\plaineqalign{
&{\lambda_j'\over|\lambda'|}\leq
{\lambda_j\over|\lambda|}\quad\hbox{for $j<n$},\qquad
{\lambda_n'\over|\lambda'|}\leq
{t'\lambda_n+n-1\over t'|\lambda|}
={\lambda_n\over|\lambda|}\bigg(1+{n-1\over t'\lambda_n}\bigg)\quad
\hbox{if $\lambda_n>0$}.\cr}
$$
These inequalities imply respectively
$$
\bigg({\lambda_j'\over|\lambda'|}\bigg)^{\lambda'_j/2}\leq
\bigg({\lambda_j\over|\lambda|}\bigg)^{(t'\lambda_j-1)/2},\quad
\bigg({\lambda_n'\over|\lambda'|}\bigg)^{\lambda'_n/2}\leq
\bigg({\lambda_n\over|\lambda|}\bigg)^{t'\lambda_n/2}
\bigg(1+{n-1\over t'\lambda_n}\bigg)^{(t'\lambda_n+n-1)/2}.
$$
In the last inequality we have $t'\lambda_n\geq{1\over|\lambda|}\lambda_n
\geq {1\over n}$ unless $\lambda'=0$. Thus, if $\lambda'\neq 0$, we get
$$
\bigg(1+{n-1\over t'\lambda_n}\bigg)^{(t'\lambda_n+n-1)/2}\leq
\exp\bigg({1\over 2}{n-1\over t'\lambda_n}(t'\lambda_n+n-1)\bigg)
\leq
\exp\Big({1\over 2}\big(n-1+n(n-1)^2\big)\Big),
$$
and by taking the product over all $j\in\{1,\ldots,n\}$ we find
$$
s_{\lambda'}\leq e^{n^3/2}\;
\prod_{j\leq n}\bigg({\lambda_j\over|\lambda|}\bigg)^{t'\lambda_j/2}
\prod_{j<n}\bigg({|\lambda|\over \lambda_j}\bigg)^{1/2}
\leq e^{n^3/2}\;(\sigma_\lambda)^{t'}\;|\lambda|^{(n-1)/2}
$$
(notice that for $\lambda_j=0$ we also have $\lambda'_j=0$, and the
corresponding factors are then equal to~$1$). Notice also that
$$
(s_\lambda)^{-1/|\lambda|}=
\prod_{j\leq n}\bigg({|\lambda|\over\lambda_j}\bigg)^{\lambda_j/2|\lambda|}
\leq\prod_{j\leq n}|\lambda|^{\lambda_j/2|\lambda|}=|\lambda|^{1/2}.
$$
For $\lambda',\lambda''\neq 0$, this implies
$$
s_{\lambda'}s_{\lambda''}\leq e^{n^3}(s_\lambda)^{t'+t''}\;|\lambda|^{n-1}
=e^{n^3}(s_\lambda)^{1-1/|\lambda|}\;|\lambda|^{n-1}
\leq e^{n^3}s_\lambda\;|\lambda|^n,\leqno(3.21)
$$
and our combinatorial factor (3.20) is less than
$e^{2n^3}|\lambda|^n|\mu|^n$. When $\lambda'=0$
or $\lambda''=0$ (say $\lambda''=0$), we have $\lambda'=\lambda-c_j$
for some $j$ and $s_{\lambda''}=1$, thus
$$
s_{\lambda'}s_{\lambda''}=s_{\lambda'}
=s_\lambda\;\bigg({|\lambda|\over|\lambda|-1}\bigg)^{(|\lambda|-1)/2}
|\lambda|^{1/2}\;
{(\lambda_j-1)^{(\lambda_j-1)/2}\over\lambda_j^{\lambda_j/2}}
\leq e^{1/2}\;s_\lambda\;|\lambda|^{1/2}
$$
and inequality (3.21) still holds. We now put all our bounds together.
For all $r'<r={}$radius of convergence of $w\mapsto\Psi_z(w)$, 
the coefficients $a'_{k,\ell,\lambda}$ satisfy
$|a'_{k,\ell,\lambda}|\leq C_0(1/r')^{|\lambda|}$ 
with $C_0=C_0(r')>0$, and for every
$\xi=\sum_\alpha\xi_\alpha\wt e_\alpha$, (3.16--3.21) imply
a bound of the form
$$
\plaineqalign{
\big|\langle Q(\xi),\xi\rangle\big|
\leq\varepsilon^{-2}\sum_{\alpha\in\bN^n}\sum_{\lambda,\mu\in\bN^n\ssm\{0\}}
C_1\;&\Big({\varepsilon\over r'}\Big)^{|\lambda|+|\mu|}
(2+|\lambda|+|\mu|+|\lambda||\mu|)\;|\lambda|^n|\mu|^n\times{}\cr
&{((|\alpha|-|\lambda|+|\mu|-1)_++n+1)(|\alpha|+n)\over
(|\alpha|-|\lambda|)_++n}|\xi_\alpha||\xi_{\alpha-\lambda+\mu}|.\cr}
$$
Here $|\lambda|+|\mu|\ge 2$, and for $\delta>0$ arbitrary, there exists
$C_2=C_2(\delta)$ such that 
$$
(2+|\lambda|+|\mu|+|\lambda||\mu|)\;|\lambda|^n|\mu|^n\leq
C_2\;(1+\delta)^{|\lambda|+|\mu|-2},
$$
thus
$$
\plaineqalign{
&\big|\langle Q(\xi),\xi\rangle\big|
\leq {C_1C_2\over r^{\prime 2}}
\sum_{\alpha\in\bN^n}\sum_{\lambda,\mu\in\bN^n\ssm\{0\}}
\Big({(1+\delta)\varepsilon\over r'}\Big)^{|\lambda|+|\mu|-2}\times{}\cr
&\kern160pt{}{((|\alpha|-|\lambda|+|\mu|)_++n)(|\alpha|+n)\over
(|\alpha|-|\lambda|)_++n}|\xi_\alpha||\xi_{\alpha-\lambda+\mu}|.\cr}
$$
Now, we split the summation with respect to $(\lambda,\mu)$ between
the two subsets $|\lambda|+|\mu|\leq(|\alpha|+n)/2$ and
$|\lambda|+|\mu|>(|\alpha|+n)/2$. We find respectively
$$
{((|\alpha|-|\lambda|+|\mu|)_++n)(|\alpha|+n)\over(|\alpha|-|\lambda|)_++n}
\leq\plaincases{
\sqrt{6}\sqrt{|\alpha|+n)((|\alpha|-|\lambda|+|\mu|)_++n)}&(first case)\cr
\noalign{\vskip4pt}
{6\over n}(|\lambda|+|\mu|)^2&(second case).\cr
}
$$
In the first case, we use the inequality
$$
\plaineqalign{
&2\sqrt{|\alpha|+n)((|\alpha|-|\lambda|+|\mu|)_++n)}\;
|\xi_\alpha||\xi_{\alpha-\lambda+\mu}|\cr
&\kern100pt{}\leq
(|\alpha|+n)|\xi_\alpha|^2+
(|\alpha|-|\lambda|+|\mu|)_++n)|\xi_{\alpha-\lambda+\mu}|^2,\cr}
$$
and in the second case we content ourselves with the simpler bound
$$
2|\xi_\alpha||\xi_{\alpha-\lambda+\mu}|\leq
|\xi_\alpha|^2+|\xi_{\alpha-\lambda+\mu}|^2.
$$
For $\varepsilon\in{}]0,r[$, the series
$$
\sum_{\lambda,\mu\in\bN^n\ssm\{0\}}
\Big({(1+\delta)\varepsilon\over r'}\Big)^{|\lambda|+|\mu|-2}\quad\hbox{and}\quad
\sum_{\lambda,\mu\in\bN^n\ssm\{0\}}
\Big({(1+\delta)\varepsilon\over r'}\Big)^{|\lambda|+|\mu|-2}(|\lambda|+|\mu|)^2
$$
can be made convergent by choosing $r'=(r+\varepsilon)/2\in{}]\varepsilon,r[$ 
and $1+\delta=\sqrt{r'/\varepsilon}$, thus there exists a positive
continuous and increasing function $\varepsilon\mapsto C(\varepsilon)$ on
$]0,r[$ such that
$$
\big|\langle Q(\xi),\xi\rangle\big|\leq C(\varepsilon)\sum_{\alpha\in\bN^n}
(|\alpha|+n)|\xi_\alpha|^2\quad\hbox{for all $\xi\in B_\varepsilon$},
$$
which is what we wanted. This bound, together with Theorem 3.10 and
the estimates from the preliminary discussion yield the following
result.

\claim 3.22. Theorem|Let $(X,\gamma)$ be a compact hermitian manifold
equipped with a real analytic metric, and let $r$ we the supremum of
the radii $r'$ of the ball bundles $\{\Vert\zeta\Vert_\gamma<r'\}$ on
which the related exponential map 
$\exph=\exph_\gamma:\{\Vert \zeta\Vert_\gamma<r'\}\subset T_X
\to X\times X$ defines a real analytic diffeomorphism
$(z,\zeta)\mapsto(z,\exph_z(\zeta))$. Then, for all
$\varepsilon<r$, the curvature tensor of the Bergman bundle
$(B_\varepsilon,h)$ satisfies an estimate
$$
\langle(\Theta_{B_\varepsilon,h}\,\xi)(v,Jv),\xi\rangle_h
=\varepsilon^{-2}\sum_{\alpha\in\bN^n}\Bigg(\bigg|
\sum_j\sqrt{\alpha_j}\;\xi_{\alpha-c_j}v_j\bigg|^2
+(1+O(\varepsilon))\sum_{j}(|\alpha|+n)\;|\xi_\alpha|^2|v_j|^2\Bigg)
$$
for every $\xi=\sum_\alpha\xi_\alpha e_\alpha\in B_{\varepsilon'}$,
$\varepsilon'>\varepsilon$, and every tangent
vector $v=\sum v_j\;\partial/\partial z_j$,
where $O(\varepsilon)=\varepsilon\;C(\varepsilon)$ for a continuous
increasing function $\varepsilon\mapsto C(\varepsilon)$ on $]0,r[$.
In particular $\Theta_{B_\varepsilon,h}$ is positive definite $($and
even coercive unbounded$\,)$ for $\varepsilon<\varepsilon_0$ small enough.
\endclaim

\claim 3.23. Remark|{\rm Under our real analyticity assumptions,
the proof makes clear that there exists an asymptotic expansion
$$
\langle(\Theta_{B_\varepsilon,h}\,\xi)(v,Jv),\xi\rangle_h
=\sum_{p=0}^{+\infty}\varepsilon^{-2+p}Q_p(z,\xi\otimes v),
$$
where
$$
Q_0(z,\xi\otimes v)=Q_0(\xi\otimes v)=
\sum_{\alpha\in\bN^n}\Bigg(\bigg|
\sum_j\sqrt{\alpha_j}\;\xi_{\alpha-c_j}v_j\bigg|^2
+\sum_{j}(|\alpha|+n)\;|\xi_\alpha|^2|v_j|^2\Bigg)
$$
corresponds to the model case $X=\bC^n$. The terms $Q_j$ can be
derived from the Taylor expansion of $\exph$ associated with the metric
$\gamma$, and they a priori depend on the coefficients of the torsion and
curvature tensor and their derivatives. In the K\"ahler case, cf.\ for
instance [Dem82, (8.5)], one has $\exph_z(\xi)=z+\xi+O(\overline z\xi^2)$
and one can check from the above calculations that $Q_1=0$.
It would be interesting to identify more precisely $Q_1$ and $Q_2$
in general. It is very likely that $Q_1$ involves the torsion form
$d\omega$ and that $Q_2$ is strongly related to the curvature tensor
of $(T_X,\omega)$.}
\endclaim

\claim 3.24. Remark|{\rm Although we have already observed that 
$B_\varepsilon$ cannot be a locally trivial holomorphic Hilbert bundle, 
as follows from Remark~3.2 and the discussion made in the introduction,
one can still endow the total space of $B_\varepsilon$ and of its Hilbert dual 
$B^\smallvee_\varepsilon$ with some sort of weird infinite dimensional complex 
space structure, for which the projections $\pi:B_\varepsilon\to X$ and
$\pi^\smallvee:B^\smallvee_\varepsilon\to X$ are holomorphic. Let us start with 
$B^\smallvee_\varepsilon$.
This space has a lot of global ``holomorphic functions'', that actually
separate all points of $B^\smallvee_\varepsilon$ except those of the
zero section. In fact, every global holomorphic function
$F\in\cB_{\varepsilon'}(X)$, $\varepsilon'>\varepsilon$,
gives rise to a function $\ell_F:B^\smallvee_\varepsilon\to\bC$ where
$\ell_F(\xi)=F_{|B_{\varepsilon,z}}\cdot \xi$ for
$\xi\in B^\smallvee_{\varepsilon,z}\subset B^\smallvee_\varepsilon$. More generally,
one can define a presheaf $\cO_{B^\smallvee_\varepsilon}$
of ``holomorphic functions'' on $B^\smallvee_\varepsilon$ as follows:
if $V\subset B^\smallvee_\varepsilon$ is an open set, we take
$\cO_{B^\smallvee_\varepsilon}(V)$ to be the closure in 
locally uniform topology in $V$ of the algebra generated by
the pull-backs $u\circ\pi^\smallvee$, $u\in\cO_X(\pi^\smallvee(V))$,
and by the functions $\ell_F$, $F\in\cB_{\varepsilon'}(\pi^\smallvee(V))$,
which are linear on the fibers
of $B^\smallvee_\varepsilon$. One then gets a genuine sheaf
$\cO_{B^\smallvee_\varepsilon}$ by sheafifying the above presheaf.
The construction of $\cO_{B_\varepsilon}$ is made by reversing
the roles of $B_\varepsilon$ and~$B^\smallvee_\varepsilon$
(the $\overline\partial$ operator of $B^\smallvee_\varepsilon$ being
the Von Neumann  adjoint of the $(1,0)$ part of the Chern connection
of $\nabla^{1,0}$ on $B_\varepsilon$, and the sheaf of ``holomorphic sections''
of $B^\smallvee_\varepsilon$ being its kernel).
}
\endclaim

\plainsection{4}{On the invariance of plurigenera for polarized
K\rlap{$\raise6.5pt\hbox{\bf.\kern-1pt.}$}ahler families}

An important unsolved problem of K\"ahler geometry is the invariance
of plurigenera for compact K\"ahler manifolds, which can be stated as follows.

\claim 4.1. Conjecture|Let $\pi:\cX\to S$ be a proper holomorphic
map defining a family of smooth compact K\"ahler manifolds over an
irreducible base~$S$. Assume that $\pi$ admits local polari\-zations,
i.e.\ every point $t_0\in S$
has a neighborhood $V$ such that $\pi^{-1}(V)$ carries a K\"ahler metric
$\omega$. Then the plurige\-nera $p_m(X_t)=h^0(X_t,mK_{X_t})$ of fibers
are independent of~$t$ for all $m\ge 0$.
\endclaim

This conjecture has been affirmatively settled by Y.T.~Siu [Siu98]
in the case of projective varieties of general type
(in which case the proof has been translated into a purely algebraic
form by Y.~Kawamata [Kaw99]), and then by [Siu02] and P{\u a}un [Pau07]
in the case of arbitrary projective varieties; remarkably, no
algebraic proof of the result is known beyond the case proved by
Kawamata. Here, we wish to study such results in the K\"ahler context.
This requires a priori substantial modifications of Siu's proof, since
the technique involves in a crucial manner the use
of an auxiliary ample line bundle. In the light of the previous sections,
a potential replacement would be to use the ``very ample'' Bergman bundles just
constructed. Conjecture 4.1 would be a consequence of the following
more technical statement.

\claim 4.2. Conjecture {\rm (generalized version of the Claudon-P\u{a}un
theorem)}|
Let $\pi: \cX \to \Delta$ be a polarized family of compact K\"ahler
manifolds over a disc~$\Delta\subset\bC$, and
let $(\cL_j,h_j)_{0\le j\le N-1}$ be $($singular$)$ hermitian
line bundles with semi-positive curvature currents $\ii\Theta_{\cL_j,h_j}\ge 0$
on~$\cX$. Assume that
\plainitem{\rm (a)} the restriction of $h_j$ to the central fiber $X_0$ is well
defined $($i.e.\ not identically~$+\infty)$.
\smallskip
\plainitem{\rm (b)} the multiplier ideal sheaf
$\cI(h_{j|X_0})$ is trivial for $1\le j\le N-1$.
\smallskip
\noindent
Then any section $\sigma$ of
$\cO(NK_{\cX}+\sum \cL_j)_{|X_0}\otimes\cI(h_{0|X_0})$ over the central
fiber $X_0$ extends into a section $\widetilde\sigma$ of
$\cO(NK_{\cX}+\sum \cL_j)$ over a certain neighborhood $\cX'=\pi^{-1}(\Delta')$
of $X_0$, where $\Delta'\subset\Delta$ is a sufficienty small disc
centered at~$0$.
\endclaim

\noindent
The invariance of plurigenera is the special case of Conjecture
4.2~when all line
bundles $\cL_j$ and their metrics $h_j$ are trivial. Since the
dimension $t\mapsto
h^0(X_t,mK_{X_t})$ is always upper semicontinuous and since Conjecture~4.2
implies the lower semicontinuity, we conclude that the dimension
must be constant along analytic discs, hence along the irreducible base $S$,
by joining any two points through a chain of analytic discs.\qed

\claim 4.3. Remark|{\rm A standard cohomological argument shows that we
can in fact take $\cX'=\cX$ in the conclusion of Conjecture 4.2, because
the direct image sheaf $\cE=\pi_*\cO(mK_{\cX}+\sum \cL_j)$
is coherent, and the restriction
$\cE\to\cE\otimes(\cO_\Delta/\gm_0\cO_\Delta)$ induces
a surjective map at the $H^0$ level on the Stein space $\Delta$, so
we can extend $\widetilde\sigma$ mod $\pi^*\gm_0$ to $\cX$.}
\endclaim

\noindent We now indicate how the technology of Bergman bundles could possibly
be used to approach the conjectures.

\claim 4.4. Lemma|Let $\cX'=\pi^{-1}(\Delta')\to\Delta'$ be the restriction
of $\pi:\cX\to\Delta$ to a disc $\Delta'\compact\Delta$ centered at~$0$,
of radius $R'<R$. For $\varepsilon\leq\varepsilon_0=\varepsilon_0(R')$
small enough, one can find a Stein open subset $\cU'_\varepsilon\subset
\cX'\times\overline\cX$, such that the projection
$\pr_1:\cU'_\varepsilon\to \cX'$ is a complex ball bundle over
$\cX'$ that is locally trivial real analytically.
\endclaim

\proof{Proof} The arguments are very similar to those of \S1, except
for the fact that $\cX$ is no longer compact, but this is not a problem
since $\cX\to\Delta$ is proper, and since we can always shrink $\Delta$
a little bit to achieve uniform bounds (would they be needed).
Let $\gamma$ be a real analytic hermitian metric on~$\cX$, and
$\exph:T_{\cX}\to\cX$ be the corresponding real analytic
and fiber-holomorphic exponential map associated with $\gamma$, as
in \S1.  The map $\exph$ is no longer everywhere defined, but if
we restrict it to the $\varepsilon$-tubular neighborhood
of the zero section in $T_{\cX'}$, we get for $\varepsilon>0$ small enough
a real analytic diffeomorphism $(z,\xi)\mapsto(z,\exph_z(\xi))$ onto a tubular
neighborhood of the diagonal of $\cX'\times\cX'$. The rest of the
proof is identical to what we did in~\S1, taking
$$
\cU'_\varepsilon=\big\{(z,w)\in\cX'\times\overline\cX\,;\;
|\logh_z(\overline w)|_\gamma<\varepsilon\big\}.\lreqno{(4.5)}{\square}
$$

In order to study Conjecture~4.2, we first state a technical 
extension theo\-rem needed for the proof, which is a special case
of the well-known and extremely powerful Ohsawa-Takegoshi theorem
[OhT87], see also [Ohs88, Ohs94], [Dem00], the general K\"ahler case
stated below being due to [Cao17].

\claim 4.6. Proposition|Let $\pi:\cZ\to\Delta$ be a smooth and proper morphism
from a $($non compact$\,)$ K\"ahler manifold $\cZ$ to a disc
$\Delta\subset\bC$ and let $(\cL,h)$ be
a  $($singular$)$ hermitian line bundle with semi-positive
curvature current $\ii\Theta_{\cL,h}\ge 0$ on~$\cZ$.
Let $\omega$ be a global K\"ahler metric on $\cZ$, and let $dV_\cZ$,
$dV_{Z_0}$ the respective induced volume elements on $\cZ$ and
$Z_0=\pi^{-1}(0)$. Assume that $h_{Z_0}$ is well defined $($i.e.\ almost
everywhere finite$)$. Then any holomorphic section $s$ of
$\cO(K_\cZ+\cL)\otimes\cI(h_{|Z_0})$ extends into a section $\wt s$ over $\cZ$
satisfying an $L^2$ estimate
$$
\int_{\cZ}\Vert \wt s\,\Vert^2_{\omega\otimes h}dV_\cZ\le C_0
\int_{Z_0}\Vert s\Vert^2_{\omega\otimes h}dV_{Z_0},
$$
where $C_0\ge 0$ is some universal constant $($depending on $\dim\cZ$ and $\diam\Delta$, but otherwise independent of $\cZ$, $\cL$, $\ldots\,)$.
\endclaim

\claim 4.7. Remark|{\rm The assumptions of Proposition 4.6 imply that $\cZ$ is
holomorphically convex and complete K\"ahler, thus, as an alternative to the 
technique used in [Cao17], the regula\-rization arguments explained in [Dem82] 
would also apply to yield the result. We leave motivated readers eventually
complete such a proof.}
\endclaim

\proof{Attempt of proof of Conjecture $4.2$} Let $p=\pr_1:\cU'_\varepsilon\to
\cX'$ be as in Lemma~4.4, and $q=\pr_2:\cU'_\varepsilon\to\overline\cX$. 
We take $\varepsilon<\varepsilon_0$ and
use on $\cZ:=\cU'_\varepsilon$ a K\"ahler metric $\omega_0$ defined
on the Stein manifold $\cU'_{\varepsilon_0}$. On can define e.g.\ $\omega_0$
as the $\ii\ddbar$ of a strictly plurisubharmonic exhaustion function
on $\cU'_{\varepsilon_0}$, but we can also take the restriction 
of $\pr_1^*\omega+\pr_2^*\smash{\overline\omega_{|\overline\cX}}$
where $\omega$ is the K\"ahler metric on the total space $\cX$,
and $\overline\omega={}{-}\,\omega$ the corresponding K\"ahler metric on
the conjugate space $\overline\cX$.
\medskip

\noindent {\it First step: construction of a sequence of extensions
on $\cZ=\cU'_\varepsilon$ via the Ohsawa-Takegoshi extension theorem.}

The strategy is to apply iteratively the special case 4.6 of the
Ohsawa-Takegoshi extension theorem on the total space of the fibration
$$
\pi'=\pi\circ p:\cZ=\cU'_\varepsilon\to\cX'\to\Delta',
$$
and to extend sections of ad hoc pull-backs $p^*\cG$
from the zero fiber $Z_0=\pi^{\prime\,-1}(0)=p^{-1}(X_0)$
to the whole of $\cZ=\cU'_\varepsilon$. We write $h_j=e^{-\varphi_j}$ in terms
of local plurisubharmonic weights, and define
inductively a sequence of line bundles $\cG_m$ by putting
$\cG_0=\cO_{\cX'}$ and
$$
\cG_m=\cG_{m-1}+K_{\cX'}+\cL_r\qquad\hbox{if $m=Nq+r$,~~$0\le r\le N-1$}.
$$
By construction we have 
$$
\plaineqalign{
\cG_m&=mK_{\cX'}+\cL_1+\cdot\cdot\cdot+\cL_m,\kern28pt\hbox{for $1\le m\le N-1$}\,,\cr
\cG_{m+N}-\cG_m=\cG_N&=NK_{\cX'}+\cL_0+\cdot\cdot\cdot+\cL_{N-1}\,,\quad
\hbox{for all $m\geq 0$.}\cr}
$$
The game is to construct inductively families of sections, say
$\{\wt f_j^{(m)}\}_{j=1,\ldots, J(m)}$, of $p^*\cG_m$ over~$\cZ$, together
with ad hoc $L^2$ estimates, in such a way that
{\plainitemindent=9mm
\smallskip
\plainitem{\rm (4.8)} for $m=0,\ldots,N-1$, $p^*\cG_m$ is generated by $L^2$
sections $\{\wt f_j^{(m)}\}_{j=1,\ldots, J(m)}\;$ on $\cU'_{\varepsilon_0}\,$;
\smallskip
\plainitem{\rm (4.9)} we have the $m$-periodicity relations
$J(m+N)=J(m)$ and $\wt f_j^{(m)}$ is an extension of
$f_j^{(m)}:=(p^*\sigma)^qf_j^{(r)}$ over $\cZ$ for $m=Nq+r$,
where $f_j^{(r)}:=\wt f^{(r)}_{j|Z_0}$, $0\le r\le N-1$.
\smallskip}
\noindent
Property (4.8) can certainly be achieved since $\cU'_{\varepsilon_0}$ is Stein,
and for $m=0$ we can take $J(0)=1$ and $\widetilde f^{(0)}_1=1$.
Now, by induction, we equip $p^*\cG_{m-1}$ with the
tautological metric
$|\xi|^2/\sum_\ell|\wt f_\ell^{(m-1)}(x)|^2$, and
$$
\widetilde\cG_m:=p^*\cG_m-K_\cZ=p^*\cG_m-(p^*K_{\cX'}+q^*K_{\overline\cX})=
p^*(\cG_{m-1}+\cL_r)-q^*K_{\overline\cX}
$$
with that metric multiplied by $p^*h_r=e^{-p^*\varphi_r}$ and a fixed smooth
metric $e^{-\psi}$ of positive curvature on
$(-q^*K_{\overline\cX})_{|\cU'_{\varepsilon_0}}$
(remember that $\cU'_{\varepsilon_0}$ is Stein!).
It is clear that these metrics have semi-positive curvature currents on $\cZ$
(by adjusting $\psi$, we could even take them to be strictly positive
if we wanted).
In this setting, we apply the Ohsawa-Takegoshi theorem to the line bundle
\hbox{$K_\cZ+\widetilde\cG_m=p^*\cG_m$}, and extend in this way
$f_j^{(m)}$ into a section $\wt f_j^{(m)}$ over~$\cZ$. By construction
the pointwise norm of that section
in $p^*\cG_{m|Z_0}$ in a local trivialization of the bundles involved is
the ratio
$$
{|f_j^{(m)}|^2\over\sum_\ell|f_\ell^{(m-1)}|^2}e^{-p^*\varphi_r-\psi},
$$
up to some fixed smooth positive factor depending only on the metric induced
by $\omega_0$ on~$K_\cZ$. However, by the induction relations, we have
$$
{\sum_j|f_j^{(m)}|^2\over\sum_\ell|f_\ell^{(m-1)}|^2}e^{-p^*\varphi_r}=
\plaincases{
\displaystyle
{\sum_j|f_j^{(r)}|^2\over\sum_\ell|f_\ell^{(r-1)}|^2}e^{-p^*\varphi_r}
&\hbox{for $m=Nq+r$, $0<r\le N-1$},\cr
\displaystyle
{\sum_j|f_j^{(0)}|^2\over\sum_\ell|f_\ell^{(N-1)}|^2}
|p^*\sigma|^2e^{-p^*\varphi_0}
&\hbox{for $m\equiv 0\mod N$, $m>0$}.\cr}
$$
Since the sections $\{f_j^{(r)}\}_{0\leq r<N}$ generate their line bundle on
$\cU_{\varepsilon_0}\supset\overline{\cU'_\varepsilon}$, the ratios
involved are positive functions without zeroes and poles, hence smooth and
bounded [possibly after shrinking a little bit the base disc $\Delta'$,
as is permitted].
On the other hand, assumption 4.2~(b) and the fact that $\sigma$ has
coefficients in the multiplier ideal sheaf $\cI(h_{0|X_0})$ tell us that
$\smash{e^{-p^*\varphi_r}}$, $1\le r<m$ and
$|p^*\sigma|^2\smash{e^{-p^*\varphi_0}}$ are locally integrable on~$Z_0$.
It follows that there is a constant $C_1=C_1(\varepsilon)\ge 0$ such that
$$
\int_{Z_0}
{\sum_j|f_j^{(m)}|^2\over\sum_\ell|f_\ell^{(m-1)}|^2}e^{-p^*\varphi_r-\psi}
dV_{\omega_0} \le C_1
$$
for all $m\ge 1$ (of course, the integral certainly involves finitely many
trivializations of the bundles involved, whereas the integrand expression is
just local in each chart). Inductively, the $L^2$ extension theorem produces
sections $\wt f_j^{(m)}$ of $p^*\cG_m$ over $\cZ$ such that
$$
\int_{\cZ}
{\sum_j|\wt f_j^{(m)}|^2\over\sum_\ell|\wt f_\ell^{(m-1)}|^2}e^{-p^*\varphi_r-\psi}
dV_{\omega_0} \le C_2=C_0C_1.
$$

\noindent
{\it Second step: applying the H\"older inequality.} Put
$k=Nq(k)+r(k)$ with $0\leq r(k)<N$, and
take $m=Nq(m)$ to be a multiple of~$N$. 
The H\"older inequality $|\int\prod_{1\leq k\leq m}u_kd\mu|\leq
\prod_{1\leq k\leq m}(\int|u_k|^md\mu)^{1/m}$ applied to the measure
$\mu=dV_{\omega_0}$ and to the product of functions
$$
\Bigg({\sum_j|\wt f_j^{(m)}|^2\over\sum_\ell|\wt f_\ell^{(0)}|^2}\Bigg)^{1/m}
e^{-{1\over N}p^*(\varphi_0+\ldots+\varphi_{N-1})-\psi}
=\prod_{1\leq k\leq m}
\Bigg({\sum_j|\wt f_j^{(k)}|^2\over\sum_\ell|\wt f_\ell^{(k-1)}|^2}
e^{-p^*\varphi_{r(k)}-\psi}\Bigg)^{1/m}
$$
in which $\sum_\ell|\wt f_\ell^{(0)}|^2=|\wt f_1^{(0)}|^2=1$ and
$\sum_j|\wt f_j^{(m)}|^2=|\wt f_1^{(m)}|^2$, implies that
$$
\int_{\cZ}\big|\wt f_1^{(m)}\big|^{2/m}
e^{-{1\over N}p^*(\varphi_0+\ldots+\varphi_{N-1})-\psi}dV_{\omega_0}\le C_2.
\leqno(4.10)
$$
As the functions
$\smash{\varphi_{r(k)}}$ and $\psi$ are locally bounded from above,
we infer from this the weaker inequality
$$
\int_{\cZ}\big|\wt f_1^{(m)}\big|^{2/m}dV_{\omega_0}\le C_3.
\leqno(4.10')
$$
The last inequality is to be understood as an inequality that holds in fact
only locally over~$\cX'$, on sets of the form $p^{-1}(V)$, where $V\compact\cX'$
are small coordinate open sets where our line bundles are trivial, so that
the section $\wt f_1^{(m)}$ of $q(m)\,p^*(NK_{\cX'}+\sum \cL_j)$ can be viewed
as a holomorphic function on $p^{-1}(V)$.\medskip

\noindent
{\it Third step: construction of singular hermitian metrics on
$NK_{\cX'}+\sum\cL_j$}. The rough idea is to extract a weak limit of the
$m$-th root occurring in (4.10), $(4.10')$, combined with an integration
on the fibers of $p:\cZ=\cU'_\varepsilon\to\cX'$, to get a singular
hermitian metric on~\hbox{$NK_{\cX'}+\sum \cL_j$}. This is the crucial
step in the proof, and the
place where the K\"ahler setup will require new arguments; especially, the
integration on fibers makes the weak limit argument much less obvious
than in the projective setup. Our (yet incomplete) attempt involves
the results of \S2, \S3 on Bergman bundles.

\claim 4.11. Proposition|Assume that the sections $\widetilde f^{(m)}_1$
have been constructed on $\cZ=\cU'_{\varepsilon}\to\cX'$,
$\varepsilon\le \varepsilon_0(R')$, and let us shrink
these sections to a smaller neighborhood $\cU'_{\rho\varepsilon}$, $\rho<1$.
Then there exists 
a subsequence $m\in M_0\subset\bN$ such that, with respect to local
trivializations of the $\cL_j$ and local holomorphic sections 
$dw=dw_1\wedge\ldots\wedge dw_{n+1}$ of~$K_{\overline\cX})$, we have a 
well defined limit
$$
\theta(z)=\lim_{m\in M_0\atop m\to+\infty}{1\over m}
\log\int_{w\in \cU'_{\rho\varepsilon,z}}
|\wt f_1^{\kern1pt(m)}(z,w)\big|^2\,i^{(n+1)^2}dw\wedge d\overline w,\quad z\in\cX'
$$
that exists almost everywhere on~$\cX'$, and $H=e^{-N\theta}$ defines a
singular hermitian metric on $p^*(NK_{\cX'}+\sum\cL_j)$ satisfying the
following estimates$\,:$
\smallskip
\plainitem{\rm(a)}$|\sigma|^2_H=|\sigma|^2e^{-N\theta}=1$ on $X_0\subset \cX'\;;$
\smallskip
\plainitem{\rm(b)} 
$\displaystyle
\int_{\cX'}e^{-\theta}e^{-{1\over N}(\varphi_0+\ldots+\varphi_{N-1})}
dV_{\omega}<\infty\;;$
\smallskip
\plainitem{\rm(c)} there are constants $C_4,C_5>0$ such that $\theta\le C_4$ and
$\displaystyle\ii\ddbar\theta\ge -{C_5\over \varepsilon^2\rho^2}
\big(C_4-\theta\big)\,\omega.$
\endclaim

\proof{Proof} First notice that the choice of the $w$ local coordinates
on $\overline\cX$ is irrelevant in the definition of $\theta$
(the $L^2$ integrals may eventually change by bounded multiplicative factors,
which get killed as $m\to+\infty$). We apply
the mean value inequality for plurisubharmonic functions, applied on
$\omega_0$-geodesic balls of $\cZ$ centered at points
$(z,w)\in\cU'_{\rho\varepsilon}$ and of 
radius ${1\over 2}(1-\rho)\varepsilon$ (say).
As $\dim\cZ=2(n+1)$, we obtain by $(4.10')$ a uniform upper bound
$$
\plainleqalignno{
\sup_{\cU'_{\rho\varepsilon,z}}|\wt f_1^{(m)}|^{2/m}
&\leq{C_6\over((1-\rho)\varepsilon)^{4(n+1)}}
\int_{\cU'_\varepsilon}|\wt f_1^{(m)}|^{2/m}\,i^{(n+1)^2}\,
dw\wedge d\overline w\cr
&\leq{C_7\over((1-\rho)\varepsilon)^{4(n+1)}},
\quad\forall z\in\cX'.&(4.12)\cr}
$$
Here our sections can be seen as functions only locally over trivializing
open sets of the line bundles in $\cX'$, but we can arrange that there
are only finitely many of these; hence the transition automorphisms
only involve bounded constants, after raising to power~$1/m$. At 
this point, we consider the Bergman 
bundle $B_\varepsilon\to\cX'$, and write locally over~$\cX'$
$$
\wt f_1^{(m)}(z,w)\;dw=\sum_{\alpha\in\bN^{n+1}}\xi_{m,\alpha}(z)\;
\wt e_\alpha(z,w)\otimes g(z)^{q(m)},\quad z\in\cX',~
w\in\cU'_{\varepsilon,z}
$$
in terms of an orthonormal frame $(\wt e_\alpha)_{\alpha\in\bN^{n+1}}$ of
$B_\varepsilon$, of the corresponding Hilbert space coefficients 
$\xi_m=(\xi_{m,\alpha})_{\alpha\in\bN^{n+1}}$ 
as defined in \S2, and of a local holomorphic generator $g$
of $\cO_\cX(NK_{\cX'}+\sum\cL_j)$. If we put
$dw=dw_1\wedge\ldots\wedge dw_{n+1}$ in local coordinates,
we get an equality
$$
\plainleqalignno{
\theta_{m,\rho}(z):\kern-2pt&={1\over m}\log\int_{w\in \cU'_{\rho\varepsilon,z}}
|\wt f_1^{(m)}(z,w)\big|^2\,i^{(n+1)^2}dw\wedge d\overline w\cr
&={1\over m}\log\Bigg(\sum_{\alpha\in\bN^{n+1}}
\rho^{2(|\alpha|+n+1)}\;|\xi_{m,\alpha}(z)|^2\Bigg),&(4.13)\cr}
$$
and by (4.12), we obtain an upper bound
$$
\theta_{m,\rho}(z)\leq {1\over m}\log
{C_8\;(\rho\varepsilon)^{2(n+1)}\;C_7^m\over((1-\rho)\varepsilon)^{4m(n+1)}}
\leq C_9+4(n+1)\,\log{1\over(1-\rho)\varepsilon}=: C_{10,\rho,\varepsilon}.
\leqno(4.14)
$$
The sum $\sum_{\alpha\in\bN^{n+1}}\rho^{2(|\alpha|+n+1)}\;|\xi_{m,\alpha}(z)|^2=
e^{m\theta_{m,\rho}(z)}$ is nothing else than the square of the norm 
of the section $\smash{\wt f^{\,(m)}_1}$, expressed with respect to
the natural hermitian metric
$\langle{\scriptstyle\bullet},{\scriptstyle\bullet}\rangle_\rho$
of the Bergman bundle $B_{\rho\varepsilon}$.
The inequalities (4.12) show that the series converges uniformly over the
whole of~$\cX'$. As $\nabla^{0,1}\xi=0$,
a~standard calculation with respect to the Bergman connection $\nabla
=\nabla^{1,0}+\nabla^{0,1}$ of $B_{\rho\varepsilon}$ implies
$$
\plainleqalignno{
\ii\ddbar\theta_{m,\rho}
&={\ii\over m\,\Vert\xi_m\Vert_\rho^2}\Bigg(\!
\langle\nabla^{1,0}\xi_m,\kern-1pt\nabla^{1,0}\xi_m\rangle_\rho-
\langle\Theta_{B_{\rho\varepsilon}}\xi_m,\xi_m\rangle_\rho-
{\langle\nabla^{1,0}\xi_m,\xi_m\rangle_\rho\wedge
\overline{\langle\nabla^{1,0}\xi_m,\xi_m\rangle_\rho}\over
\Vert\xi_m\Vert^2_\rho}\,\Bigg)\kern-3pt\cr
&\geq-{1\over m}\;{\langle\ii\Theta_{B_{\rho\varepsilon}}\xi_m,\xi_m\rangle_\rho\over
\Vert\xi_m\Vert^2_\rho}&(4.15)\cr}
$$
by the Cauchy-Schwarz inequality. On the other hand, 
as the orthonormal coordinates expressed in $B_{\rho\varepsilon}$ are
the $(\rho^{|\alpha|+n+1}\xi_{m,\alpha})$, the curvature bound
obtained in \S2 yields
$$
\langle\ii\Theta_{B_{\rho\varepsilon}}\xi_m,\xi_m\rangle_\rho
\leq (2+O(\rho\varepsilon))(\rho\varepsilon)^{-2}
\sum_{\alpha\in\bN^{n+1}}(|\alpha|+n+1)\,\rho^{2(|\alpha|+n+1)}\;|\xi_{m,\alpha}|^2\;
\omega.
$$
The last two inequalities imply the fundamental estimate
$$
\plainleqalignno{
\ii\ddbar\theta_{m,\rho}
&\geq-{(2+O(\rho\varepsilon))(\rho\varepsilon)^{-2}\over m}~
{\displaystyle
\sum_{\alpha\in\bN^{n+1}}(|\alpha|+n+1)\,\rho^{2(|\alpha|+n+1)}\;|\xi_{m,\alpha}|^2\over
\displaystyle
\sum_{\alpha\in\bN^{n+1}}\rho^{2(|\alpha|+n+1)}\;|\xi_{m,\alpha}|^2}\;\omega
&(4.16)\cr
&\geq-{1+O(\rho\varepsilon)\over \varepsilon^2\rho}
\bigg({\partial\over\partial\rho}\theta_{m,\rho}\bigg)\;\omega.&(4.16')\cr}
$$
From its definition, we see that $\theta_{m,\rho}$ is a convex function
of $\log\rho$. Therefore, for $\rho\leq \rho_1<1$, we have
$$
\rho\;{\partial\over\partial\rho}\theta_{m,\rho}\leq
{\theta_{m,\rho_1}-\theta_{m,\rho}\over\log\rho_1-\log\rho}
\leq{C_{9,\rho_1,\varepsilon}-\theta_{m,\rho}\over \log\rho_1},
$$
by (4.14), and $(4.16')$ implies
$$
\ii\ddbar\theta_{m,\rho}\geq-{C_{11}\over \varepsilon^2\rho^2}
\big(C_{10,\rho_1,\varepsilon}-\theta_{m,\rho}\big)\,\omega.
$$
A straightforward calculation yields
$$
-\ii\ddbar\log(C_{10,\rho_1,\varepsilon}+1-\theta_{m,\rho})
\geq-{C_{11}\over \varepsilon^2\rho^2}\,\omega,
$$
hence the functions $u_m=-\log(C_{10,\rho_1,\varepsilon}+1-\theta_{m,\rho})\le 0$
have Hessian forms that are uniformly bounded
from below. Also, by construction (cf.\ 4.9), $\theta_{m,\rho}$ converges to
${1\over N}\log|\sigma|$ on $X$. Standard results of pluripotential
theory imply that we can find a subsequence of $(u_m)$ that converges in $L^p$
topology (for every $p\in{}[1,+\infty[$) and pointwise almost everywhere.
Therefore we can find a limit $\theta_{m,\rho}\to \theta$
satisfying the Hessian estimates
$$
\ii\ddbar \theta\geq-{C_{11}\over \varepsilon^2\rho^2}
\big(C_{10,\rho_1,\varepsilon}-\theta\big)\,\omega,\qquad
-\ii\ddbar\log(C_{10,\rho_1,\varepsilon}+1-\theta)
\geq-{C_{11}\over\varepsilon^2\rho^2}\,\omega
$$
Proposition 4.11 is proved, as estimate (b) follows from (4.10).
\medskip

\noindent
{\it Fourth step: applying Ohsawa-Takegoshi once again with the singular
hermitian metric produced in the third step.}
Assume that we can replace estimate 4.11~(c) by the stronger fact
that the curvature form of $H=e^{-N\theta}$ is positive in the sense of
currents, i.e.
$$
-\ii\ddbar\log H=N\ii\ddbar\theta\ge 0.\leqno(4.17)
$$
This means that $NK_{\cX'}+\sum \cL_j$
possesses a hermitian metric $H$ such that $\Vert\sigma\Vert_H\le 1$
on $X_0$ and $\Theta_H\ge 0$ on~$\cX'$. 
In~order to conclude, we proceed as Siu and P\u{a}un, and equip the bundle
$$
\cE=(N-1)K_{\cX'}+\sum \cL_j
$$
with the metric $\eta=H^{1-1/N}\prod h_j^{1/N}$, and
$NK_{\cX'}+\sum \cL_j=K_{\cX'}+\cE$ with the metric $\omega\otimes\eta$.
It is important here that $\cX$ possesses a global K\"ahler polarization
$\omega$, otherwise the required estimates would not be valid.
Clearly $\eta$ has a semi-positive curvature current on $\cX'$ and in a
local trivialization we have
$$
\Vert\sigma\Vert^2_{\omega\otimes\eta}\le
C|\sigma|^2\exp\Big(-(N-1)\theta-{1\over N}\sum \varphi_j\Big)
\le C \Big(|\sigma|^2\prod e^{-\varphi_j}\Big)^{1/N}
$$
on $X_0$. Since $|\sigma|^2e^{-\varphi_0}$ and $e^{-\varphi_r}$, $r>0$ are
all locally integrable, we see that $\Vert\sigma\Vert^2_{\omega\otimes\eta}$
is also locally integrable on $X_0$ by the H\"older inequality. A new (and
final) application of the
$L^2$ extension theorem to the hermitian line bundle $(\cE,\eta)$ implies
that $\sigma$ can be
extended to~$\cX'$. Conjecture 4.2 would then be proved.
\endproof

\noindent
{\it Fifth step: final discussion.} Unfortunately, estimate (4.17) will
a priori hold only in the case where $\varepsilon$ can be taken arbitrarily
large (in the sense that the exponential map is at least everywhere an
immersion -- one can then argue on the ``unfolded neighborhood''
$\smash{\widetilde U_\varepsilon}$ diffeomorphic to the $\varepsilon$-tubular
neihborhood of the $0$ section in $T_X$, equipped with the complex
structure obtained by pulling back the complex structure of
$X\times\ol X$ via $\exph$. This condition is met e.g.\ when $X$ is a
complex torus or a ball quotient.
However, it is doubtful that all compact K\"ahler manifolds with
$K_X$ pseudo-effective satisfy this property. The main issue is that
the unboundedness of $\Theta_{B_\varepsilon,h}$ does not a priori imply
that the right hand side of (4.15) converges weakly to $0$, while this is
obviously true in the algebraic situation where we use instead a given ample
line bundle~$A$ on $\cX$. One possible way to circumvent this
difficulty is to observe that the term
$\langle\ii\Theta_{B_{\rho\varepsilon}}\xi_m,\xi_m\rangle_\rho$ is
controlled by $\Vert\xi_m\Vert_\rho \Vert\xi_m\Vert'_\rho$ where
$$
\Vert\xi_m\Vert_\rho^{\prime\,2}:=
\sum_{\alpha\in\bN^{n+1}}(|\alpha|+n+1)^2\,\rho^{2(|\alpha|+n+1)}\;
|\xi_{m,\alpha}|^2\sim
\int_{w\in \cU'_{\rho\varepsilon,z}}|\wt f_1^{(m)}(z,w)\big|^2+
|D_w\wt f_1^{(m)}(z,w)\big|^2,
$$
and it would be sufficient to find extensions $\wt f_1^{(m)}$
satisfying the additional estimate
$$
\int_{w\in \cU'_{\rho\varepsilon,z}}|D_w\wt f_1^{(m)}(z,w)\big|^2\leq
K_m\int_{w\in \cU'_{\rho\varepsilon,z}}|\wt f_1^{(m)}(z,w)\big|^2
\leqno(4.18)
$$
where $K_m$ grows subquadratically, i.e.\ ${1\over m^2}K_m\to 0$.
Getting such an estimate, e.g.\ a bound
$K_m=O(m)$  in the general situation, does not appear to be
completely implausible, since the main inductive step consists of
extending a section multiplied by $p^*\sigma(z,w)=\sigma(z)$, which is therefore
independent of~$w$ on $\cX'_0$. In this process, one might hope to obtain an
appropriate $L^2$ extension theorem taking
care of ``vertical derivatives'' with respect to a given morphism
$\cY\to\cX\to\Delta$ (namely, $\cU'_\varepsilon\to\cX'\to\Delta'$ in this
circumstance). We will
try to investigate these questions in the near future.

\vskip30pt

\centerline{\twelvebf References}
\vskip8pt\parskip=3.5pt plus 1 pt minus 1pt

\Bibitem[Ber09]&Berndtsson, B.:& Curvature of vector bundles associated
  to holomorphic fibrations;& Ann.\ of Math.\ (2) {\bf 169} (2009)
  531--560&

\Bibitem[Cao17]&Cao, J.:& Ohsawa-Takegoshi Extension Theorem for Compact
  K\"ahler Manifolds and Applications;& Springer INdAM Series, 19--38&

\Bibitem[Dem82]&Demailly, J.-P.:& Estimations $L^2$ pour l'op\'erateur $\dbar$
  d'un fibr\'e vectoriel holomorphe semi-positif au dessus d'une vari\'et\'e
  k\"ahl\'erienne compl\`ete;& Ann.\ Sci.\ Ec.\ Norm.\ Sup.\ {\bf 15}
  (1982) 457--511&
  
\Bibitem[Dem94]&Demailly, J.-P.:& Regularization of closed positive
  currents of type $(1,1)$ by the flow of a Chern connection;& Actes du
  Colloque en l'honneur de P.~Dolbeault (Juin 1992), \'edit\'e
  par H.~Skoda et J.M.~Tr\'epreau, Aspects of Mathematics, Vol.~E~26,
  Vieweg (1994), 105--126&

\Bibitem[Dem00]&Demailly, J.-P.:&On the Ohsawa-Takegoshi-Manivel $L^2$
  extension theorem;& Proceedings of the Conference in honour of the
  85th birthday of Pierre Lelong, Paris, September 1997, \'ed.\ P.~Dolbeault,
  Progress in Mathematics, Birkh\"auser, Vol.~{\bf 188} (2000) 47--82&

\Bibitem[Gra58]&Grauert, H.:& On Levi's Problem and the imbedding of
  real-analytic manifolds;& Annals of Math.\ {\bf 68} (1958) 460--472&

\Bibitem[Gro55]&Grothendieck, A.:& Produits tensoriels topologiques et
  espaces nucl\'eaires;& Providence, Memoir Amer.\ Math.\ Soc.\
  nr.~16, 1955&

\Bibitem[Hor66]&H\"ormander, L.:& An introduction to Complex Analysis in
  several variables;& 1st edition, Elsevier Science Pub., New York, 1966,
  3rd revised edition, North-Holland Math.\ library, Vol 7, Amsterdam (1990)&

\Bibitem[Kaw99]&Kawamata, Y.:& Deformation of canonical singularities;&
  J.\ Amer.\ Math.\ Soc.\ {\bf 12} (1999), 85--92&

\Bibitem[Kod54]&Kodaira, K.:& On K\"ahler varieties of restricted type
  $($an intrinsic characterization of algebraic varieties$\,)$;& 
  Annals of Math.\ {\bf 60} (1954) 28--48&

\Bibitem[Lae02]&Laeng, L.:& Estimations spectrales asymptotiques en
  g\'eom\'etrie hermitienne;& PhD thesis defended at Universit\'e
  Joseph-Fourier Grenoble I, October 2002&

\Bibitem[LeS14]&Lempert, L., Sz\H{o}ke, R.:& Direct Images, fields of
  Hilbert spaces, and geometric quantization;& Comm.\ in Math.\ Phys.\
  {\bf 327} (2014) 49-–99&

\Bibitem[Mal58]&Malgrange, B.:& Lectures on the theory of functions of
  several complex variables;& Tata Institute of Fundamental Research,
  Bombay, 1958, reprinted by the Tata Institute, 1962&

\Bibitem[Ohs88]&Ohsawa, T.:& On the extension of $L^2$ holomorphic
  functions, II;& Publ.\ RIMS, Kyoto Univ.\ {\bf 24} (1988) 265--275&

\Bibitem[Ohs94]&Ohsawa, T.:&On the extension of $L^2$ holomorphic functions,
  IV$\,$: A new density concept;& Mabuchi, T.\ (ed.) et al.,
  Geometry and analysis on complex manifolds. Festschrift for Professor
  S.~Kobayashi's 60th birthday. Singapore: World Scientific,
  (1994) 157--170&

\Bibitem[OhT87]&Ohsawa, T., Takegoshi, K.:& On the extension
  of $L^2$ holomorphic functions;& Math.\ Zeitschrift {\bf 195} (1987)
  197--204&

\Bibitem[Pau07]&P\u{a}un, M.;& Siu's invariance of plurigenera: a one-tower
  proof;& J.~Differential Geom.\ {\bf 76} (2007) 485--493&

\Bibitem[Pop13]&Popovici, D.:& Transcendental K\"ahler cohomology classes;&
  Publ.\ Res.\ Inst.\ Math.\ Sci.\ {\bf 49} (2013) 313--360&

\Bibitem[Siu98]&Siu, Y.T.;& Invariance of plurigenera;&
  Invent.\ Math.\ {\bf 134} (1994) 631--639&
  
\Bibitem[Siu02]&Siu, Y.-T;& Extension of twisted pluricanonical sections 
  with plurisubharmonic weight and invariance of semipositively twisted 
  plurigenera for manifolds not necessarily of general type;& Complex 
  Geometry (G\"ottingen, 2000), Springer, Berlin, 2002, 223--277&

\Bibitem[ShZ99]&Shiffman, B., Zelditch, S.:& Distribution of zeros of
  random and quantum chaotic sections of positive line bundles;&
  Commun.\ Math.\ Phys.\ {\bf 200} (1999) 661--683&

\Bibitem[Wan17]&Wang, Xu:& A curvature formula associated to a family of
  pseudoconvex domains;& Ann.\ Inst.\ Fourier {\bf 67} (2017) 269--313&

\Bibitem[Web89]&Webster, S.M.:& A new proof of the Newlander-Nirenberg
  theorem;& Math.\ Zeit.\ {\bf 201} (1989) 303--316&

\Bibitem[ZeZ18]&Zelditch, S., Zhou, Peng:& Pointwise Weyl law for partial
Bergman kernels;& Algebraic and analytic microlocal analysis,
589-634, Springer Proc. Math. Stat., 269, Springer, Cham, 2018&

\bigskip\noindent
Jean-Pierre Demailly\\
Universit\'e Grenoble Alpes, Institut Fourier\\
100 rue des Maths, 38610 Gi\`eres, France\\
e-mail: {\tt jean-pierre.demailly@univ-grenoble-alpes.fr}

\medskip\noindent
(revised on April 5, 2021, printed on \today, \timeofday)
\end{document}